\documentclass[english,12pt,leqno]{article}
\usepackage{tikz}
\usetikzlibrary{decorations.markings, patterns}
\usetikzlibrary{calc}

\usepackage[T2A]{fontenc}

\usepackage[latin1, russian, utf8]{inputenc}

\usepackage{amsxtra}
\usepackage{amsmath}
\usepackage{amssymb}
\usepackage{amsfonts}
\usepackage[all]{xy}
\usepackage{mathrsfs}
\usepackage{amsthm}
\usepackage{float}
 \usepackage{caption}
 \usepackage{epigraph}
 \usepackage{wrapfig}
 \usepackage{mathabx}

\newtheorem{thm}[equation]{Theorem}

\newtheorem{cor}[equation]{Corollary}
\newtheorem{lem}[equation]{Lemma}
\newtheorem{prop}[equation]{Proposition}

\newtheoremstyle{example}{\topsep}{\topsep}%
     {}
     {}
     {\bfseries}
     {.}
     {2pt}
     {\thmname{#1}\thmnumber{ #2}\thmnote{ #3}}

   \theoremstyle{example}
   
   \newtheorem{Defi}[equation]{Definition}
   
   \newtheorem{rem}[equation]{Remark}
   \newtheorem{rems}[equation]{Remarks}
   
   \newtheorem{exas}[equation]{Examples}
   \newtheorem{ex}[equation]{Example}

\newtheorem{conj}[equation]{Conjecture}

 \newtheorem{secmic}[equation]{Second Microlocalization Problem}

\newtheoremstyle{example}{\topsep}{\topsep}%
     {}
     {}
     {\bfseries}
     {.}
     {2pt}
     {\thmname{#1}\thmnumber{ #2}\thmnote{ #3}}

  \numberwithin{equation}{section}

\def\eps{{\varepsilon}}

\def\CC{\mathbb{C}}

\def\RR{\mathbb{R}}

\def\HH{\mathbb{H}}

\def\hen{\mathfrak{h}}
\def\len{\mathfrak{l}}

\def\sen{\mathfrak{s}}

\def\Kc{\mathcal{K}}
\def\Dc{\mathcal{D}}
\def\Ec{\mathcal{E}}
\def\Fc{\mathcal{F}}
\def\Gc{\mathcal{G}}

\def\Lc{\mathcal{L}}
\def\Mc{\mathcal{M}}

\def\Hc{\mathcal{H}}

\def\Pc{\mathcal{P}}
\def\Qc{\mathcal{Q}}

\def\Sc{\mathcal{S}}

\def\Uc{\mathcal{U}}
\def\Vc{\mathcal{V}}

\def\<{\langle}
\def\>{\rangle}
\def\1{{\bf 1}}

 \def\Ap{{ {}^\prime \hskip -.1cm A}}

\def\be{\begin{equation}}
\def\ee{\end{equation}}

\def\bef{\begin{figure}[H]\centering}
\def\enf{\end{figure}}
 
\def\bic{{\on{bicon}}}
\def\btp{\begin{tikzpicture}}
\def\etp{\end{tikzpicture}}

\def\codim{\on{codim}}
\def\Coker{\on{Coker}}
\def\con{{\on{con}}}

 \def\emb{{\on{emb}}}

 \def\Fl{{\on{Fl}}}
\def\FS{{\on{FS}}}

\def\Hom{\on{Hom}}

\def\Hyp{{\on{Hyp}}}

\def\Id{{\on{Id}}}

\def\k {\mathbf k}
\def\Ker{\on{Ker}}

\def\Lin{\on{Lin}}
\def\lla{\longleftarrow}
\def\lra{\longrightarrow}

\def\ol{\overline}
\def\on{\operatorname}
 \def\orr{\on{or}}
\def\OR{\on{or}}

\def\Perv{{\on{Perv}}}
\def\Pol{{\on{Pol}}} 
\def\pt{{\on{pt}}}

\def\Rep{\on{Rep}}

\def\Sh{\on{Sh}}

 \def\SS{\on{SS}}

 \def\Tot{{\on{Tot}}}

\def\ul{\underline}

\def\Vect{\on{Vect}}

\def\wt{\widetilde}

\title{Fourier-Sato transform on hyperplane arrangements}

\author{ Michael Finkelberg, Mikhail Kapranov, Vadim Schechtman }

\begin{document}


 \maketitle

 \centerline{\em To our friends and teachers Sasha Beilinson and Vitya Ginzburg}

 \addtocounter{section}{-1}
 \tableofcontents
 
 \vfill\eject
 
 \section{Introduction}
 
\paragraph{A. Setup and goals.}  The theory of perverse sheaves can be said to provide an interpolation between homology and cohomology (or  to mix them in a self-dual way). Since homology, sheaf-theoretically,
 can be understood as cohomology with compact support, interesting operations on
 perverse sheaves usually combine 
  the functors of the types  $f_!$ and $f_*$ or, dually,
  the functors of the types $f^!$ and $f^*$ in the classical
 formalism of Grothendieck. 
 
 \vskip .2cm
 
 An important context when this point of view can be pushed quite far, is that of perverse sheaves
$\Fc$  on a complex affine space $\CC^n$   smooth with respect to the stratification given
 by an arrangement  $\Hc$ of hyperplanes with real equations \cite{KS}. Denoting by
 $i_\RR: \RR^n\hookrightarrow \CC^n$ the embedding, we associate to such an   $\Fc$
 its {\em hyperbolic stalks}
 \[
 E_A(\Fc) \,=\, R\Gamma(A, i_A^* i_\RR^! \Fc).
 \]
 Here $i_A: A\hookrightarrow \RR^n$ is  the embedding of a face (stratum)
  of the real arrangement.
 It is remarkable that the $E_A(\Fc)$ reduce to single vector spaces, not complexes
 (while the ordinary stalks of $\Fc$ are of course complexes, $\Fc$ being a complex of sheaves). 
 This  type of phenomena was originally observed by  T. Braden
  in the context of varieties with a  $\CC^*$-action \cite{Br}. 
  
  \vskip .2cm
  
  It was shown in \cite{KS} that the vector spaces $E_A(\Fc)$ together with natural
  linear maps $\gamma_{AB}, \delta_{BA}$ (``generalization and specialization'')
  connecting them, determine the
  perverse sheaf $\Fc$ uniquely. Moreover, the category $\Perv(\CC^n,\Hc)$
  of perverse sheaves of the above type is equivalent to the category $\Hyp(\Hc)$
  formed by linear algebra data $(E_A, \gamma_{AB}, \delta_{BA})$
  satisfying an explicit set of conditions. We call such linear algebra data
  {\em hyperbolic sheaves}, see
  \S \ref{sec:real-com}D. 
  
  \vskip .2cm
  
  The goal of this paper is to develop the beginnings of  a ``hyperbolic calculus'',
  describing the effect of several standard operations on perverse sheaves
  directly in terms of hyperbolic sheaves. These operations include
  forming vanishing cycles,
  specialization and  Fourier-Sato transform. 
  To illustrate the importance of  such questions  recall  \cite{BFS}
  that the the weight components of
  the highest weight modules (e.g.\ Verma, or  their irreducible quotients) over
   quantized Kac-Moody algebras have interpretation as the spaces of
  vanishing cycles $\Phi_f(\Fc)$ for  appropriate $\Fc\in\Perv(\CC^n, \Hc)$ and $f$.
   In this case $\Hc$ is a  so-called discriminantal arrangement,
   $\Fc$  is  an extension  of a 1-dimensional local system
  on the generic stratum, and $f$ is  a linear
  function. The monodromy of Fourier-Sato transforms of these sheaves is related
  to the action of Lusztig symmetries on the corresponding 
  representations~\cite{FS}.

  \paragraph{B. Pattern of the results.} 
    To identify  the effect of each operation on perverse sheaves above, 
  we produce a new hyperbolic sheaf 
  out of a given one.  Our constructions and results fall into the following pattern.

  \begin{itemize}
  \item[(1)]  Each vector space of the new hyperbolic sheaf
  is identified with the $0$th cohomology space of 
  an otherwise acyclic complex formed by some of the vector spaces $E_A\otimes \OR_A$
  (here $\OR_A$ is the orientation space),
  with the differential formed out of either the $\gamma_{AB}$ or  the $\delta_{BA}$.
  So there are two versions  of the answer: the $\gamma$-answer and the
  $\delta$-answer, in each case.
  
  \item[(2)] The complexes in (1) are subquotients of the two fundamental complexes
  (Proposition \ref{f-Rgamma-c}) calculating $R\Gamma_c(\CC^n,\Fc)$ and $R\Gamma(\CC^n, \Fc)$. 
  These complexes  are sums over all the faces $A$ of the spaces $E_A\otimes\OR_A$ and 
   their differentials are formed out of the $\gamma_{AB}$ and $\delta_{AB}$
  respectively. The $R\Gamma_c(\CC^n,\Fc)$ and $R\Gamma(\CC^n, \Fc)$ typically
    have more than one nonzero cohomology, but the
  subquotients we  take turn out to be acyclic outside degree $0$. 
  
  \item[(3)] The choice of subquotient is obtained by taking not all but some summands
  $E_A\otimes\OR_A$. The selection rule, depending on the problem, reflects the geometry
  of the problem in some rough (``tropical'') way. 
  
  \item[(4)] In each case there is also a companion real statement, about   complexes
  of sheaves
  on $\RR^n$ constructible w.r.t. the stratification by the faces. This real statement
 is proved first, and the statement  for perverse sheaves is deduced from it. 
  \end{itemize}
  
  \paragraph{C. Structure of the paper.} 
  In \S 1 we recall the basics of the description of $\Perv(\CC^n, \Hc)$ by hyperbolic sheaves. 
  
  \vskip .2cm
  
  \S 2 is devoted to the calculation of the space of vanishing cycles $\Phi_f(\Fc)$  in terms of hyperbolic sheaves.
  Here $f: \CC^n\to\CC$ is a linear function with real coefficients. The selection rule for subquotients
  of $R\Gamma_c(\CC^n,\Fc)$ and $R\Gamma(\CC^n, \Fc)$ consists in taking all faces $B\subset\RR^n$
  on which $f\geq 0$. 
  
  \vskip .2cm
  
  \S 3 describes the specialization of $\Fc\in\Perv(\CC^n, \Hc)$ along a $\CC$-vector subspace $L_\CC\subset \CC^n$
  with real equations. This is a perverse sheaf $\nu_L(\Fc)$ on the normal bundle $T_L\CC^n$ which is itself a vector space.
  In this case we have the real subspace $L_\RR$, and the { product arrangement} $\nu_L(\Hc)$
  in $T_{L_\RR}\RR^n$. We further have the specialization at the level of faces which is a monotone map
  of posets 
  \[
  \nu: \bigl\{ \text{faces of } \Hc\bigr\} \to \bigl\{ \text{faces of }
\nu_L(\Hc)\bigr\}.
\]   
The selection rule   for subquotients
  of $R\Gamma_c(\CC^n,\Fc)$ and $R\Gamma(\CC^n, \Fc)$  consists in taking all  faces $A$ with
$\nu(A)=B$ being a fixed face $B$ of $\nu_L(\Hc)$. This produces  complexes calculating
the   hyperbolic stalk of $\nu_L(\Fc)$ at  $B$. 
  
We also give a description of the specialization for constructible sheaves of $\RR^n$ as the direct image under an appropriate
cellular map $q: \RR^n\to  T_{L_\RR} \RR^n$. This allows us to  identify (in our particular case)
different  possible (and, in general, non-equivalent) definitions of the  bispecialization functor
 \cite{schapira-takeuchi}
   \cite{takeuchi}
for a flag of subspaces $N\subset M\subset V$. 

\vskip .2cm

In \S 4 we give a similar description of the Fourier-Sato transform $\FS(\Fc)$ which is a perverse sheaf on the dual space $(\CC^n)^*$.
It is smooth with respect to an appropriate arrangement $\Hc^\vee$. Each face $A^\vee$
on $\Hc^\vee$ gives a natural strictly convex cone $V(A^\vee)\subset\RR^n$. The selection rule
for subquotients
  of $R\Gamma_c(\CC^n,\Fc)$ and $R\Gamma(\CC^n, \Fc)$ 
consists in taking all faces $B\subset V(A^\vee)$ for a fixed $A^\vee$. 
 This produces  complexes
calculating the hyperbolic stalk of $\FS(\Fc)$ at $A^\vee$.

\vskip .2cm

Combining the descriptions of the specialization and of the Fourier-Sato transform at the level of hyperbolic sheaves,
one obtains a description of the microlocalization  $\mu_L(\Fc)$ along a linear subspace with real equations.
The  final \S 5 is dedicated to comparison, in our linear case,
 of several possible definitions of the second microlocalization of Kashiwara and Laurent,
 see  \cite{laurent}  \cite{schapira-takeuchi}
   \cite{takeuchi}.

\paragraph{Acknowledgements.} We would like to thank
 P. Schapira for remarks on a preliminary draft of the paper
and for comminicating to us a  proof of Theorem \ref{refo:micro-2}. 
We are also grateful to Peng Zhou for several corrections. 

The research of M.F. was supported by the grant RSF 19-11-00056.
 
The research of
 M.K. was supported  by the World Premier International Research Center Initiative (WPI Initiative), MEXT, Japan.
 
 V. S. thanks Kavli IPMU   for  support of a visit during  the preparation of this paper.

 \vfill\eject
 
 \section{ Real and complex data associated to  perverse sheaves}\label{sec:real-com}

\paragraph{A. The real setup.}  
Let $V_\RR = \RR^n$ be a finite-dimensional vector space over $\RR$ and $\Hc$ be a finite
 central arrangement of hyperplanes in $V_\RR$.  
 We denote by $\Sc_\RR = \Sc_{\RR, \Hc}$ the poset of {\em faces}
 of $\Hc$, see, e.g.,  \cite{KS}, \S 2A.  Faces form a real stratification of $V_\RR$ into (a disjoint union of) 
 locally closed polyhedral cones. 
 The order $\leq$ on $\Sc_\RR$ is by inclusion of closures: $A\leq B$ means $A\subset \ol B$. 
 For an integer $p\geq 0$ we  use the notation $A<_p B$ to signify that $A\leq B$ and $\dim(B)=\dim(A)+p$, in particular $A<_0 B$  means $A=B$.
 We denote by $i_A: A\to V_\RR$ the embedding of a face $A$.

 \vskip .2cm

  Let $\k$ be a field and $\Vect_\k$ be the category of finite-dimensional $\k$-vector spaces.
  For any poset $S$ we denote by $\Rep(S)$ 
    the abelian category of {\em representations of $S$ over $\k$},
   i..e, of covariant functors from $S$ (considered as a category) to $\Vect_\k$. 
   By $D^b(\Rep(S))$ we denote the bounded derived category of $\Rep(S)$.

    For a topological space $X$ we denote by $\Sh_X$ the category of sheaves of
 $\k$-vector spaces on $X$ and by $D^b(X)$ the derived category of $\Sh_X$. 
   
   \vskip .2cm
   
    We denote by $\Sh(V_\RR, \Sc_\RR)$ the abelian category formed by sheaves of
 $\k$-vector spaces on $V_\RR$ which are constructible  with respect to the stratification
 $\Sc_\RR$.
 Let also $D^b(V_\RR, \Sc_\RR)$ be the full subcategory in the bounded derived category
 of sheaves of $\k$-vector spaces on $V_\RR$ formed by complexes with all cohomology
 sheaves lying in  $\Sh(V_\RR, \Sc_\RR)$. For $\Gc\in D^b(V_\RR, \Sc_\RR)$ and a face $A$ we denote
 \be\label{eq:stalks}
 \Gc_A \,=\, R\Gamma(A, \Gc) \,:= \, R\Gamma(A, i_A^*\Gc) \,\in \, D^b(\Vect_\k)
 \ee
 the stalk of $\Gc$ at $A$. Thus $\Gc_A$ is a complex which is a single vector space, if $\Gc$
 is a single sheaf. The following is well known.
 
 \begin{prop}\label{prop:constr-gener}
 (a) We have an equivalence of categories
 \[
 \Sh(V_\RR, \Sc_\RR) \lra \Rep(\Sc_\RR), \quad \Gc\mapsto \bigl( \Gc_A,\gamma_{AB}: \Gc_A\to \Gc_{B},
 \, A\leq B
 \bigr). 
 \]
Here  $\gamma_{AB}$ is the {\em generalization map}.

\vskip .2cm

(b) The natural functor $D^b(\Sh(V_\RR, \Sc_\RR))\to D^b(V_\RR, \Sc_\RR)$ is an equivalence. In particular:

\vskip .2cm

(c) We have an equivalence of categories $D^b(\Sh(V_\RR, \Sc_\RR))\to D^b(\Rep(\Sc_\RR))$. \qed
 \end{prop}
 
 In view of (b), we can interpret the equivalence in (c) as sending a complex of sheaves $\Gc$ to the
 collection of complexes of vector spaces $\Gc_A$ defined by \eqref{eq:stalks} 
  and generalization maps (morphisms of complexes)
 $\gamma_{AB}$ connecting them.

  \vskip.2cm
  
  By a {\em cell} we mean a topological space $B$ homeomorphic to $\RR^d$ for some $d$. 
  For a cell $B$ we denote by  $\OR_B = H^{\dim(B)}_c(B, \k)$ the 1-dimensional
{\em orientation vector $\k$-space} of $B$. For two cells $B, C$ we set
$\OR_{B/C} = \OR_C\otimes\OR_B^*$ and call it the {\em relative orientation space}
of $C$ and $B$. 

\vskip .2cm

In particular, 
  any face $B\in\Sc_\RR$  is a cell and so we have the space $\OR_B$. 
   When $B,C$ are two  faces such that $B<_1 C$,
we have a canonical  ``incidence isomorphism''
\[
\eps_{BC}: \OR_B\to\OR_C.
\]
It can be seen as a canonical trivialization of $\OR_{C/B}$. 
If $B<_1  C_1, C_2 <_1 D$ is a square of codimension $1$ inclusion of faces, then the diagram
\be\label{eq:or-anti}
\xymatrix{
\OR_B 
\ar[d]_{\eps_{B, C_2}} \ar[r]^{\eps_{B, C_1}}& \OR_{C_1}\ar[d]^{\eps_{C_1, D}}
\\
\OR_{C_2} \ar[r]_{\eps_{C_2, D}} & \OR_D
}
\ee
is anti-commutative.

 \vskip .2cm
 
  Let $j_A: A\to V_\RR$ be the embedding of a face $A$. If $A<_1 A'$
  are two faces of $\Hc$, we have a canonical morphism 
  $\xi_{AA'}: j_{A!} \ul \k_A \lra j_{A' !} \ul \k_{A'}[1]$ in $D^b(V_\RR, \Sc_\RR)$.
   Viewed as an element of $\on{Ext}^1 (j_{A!} \ul \k_A, j_{A'!} \ul \k_{A'})$,
   it represents the extension given by the subsheaf in $(j_{A' })_* \,  \ul \k_{A'}$
   formed by sections which vanish on all codimension 1 faces of $A'$ except $A$. 
   The morphisms $\xi_{AA'}$ anticommute in squares of codimension 1 embeddings,
   just like the morphisms $\eps_{AA'}$ in \eqref{eq:or-anti}.

 \begin{prop}\label{prop:LA-res}
 For $\Gc\in D^b(V_\RR, \Sc_\RR)$, the following are equivalent;
 
 \vskip .2cm
 
 (i)  $\Gc$ corresponds to  the data
  $(\Gc_A,\gamma_{AB})$.
  
  \vskip .2cm
  
  (ii) We have a resolution of $\Gc$ (a complex over $D^b(Sh_{V_\RR})$ with total object $\Gc$) 
  of the form
  \[
  \bigoplus_{\dim(A)=0} j_{A!} (\ul{\Gc_A}_A) \buildrel\gamma\otimes \xi\over\lra 
    \bigoplus_{\dim(A)=1} j_{A!} (\ul{\Gc_A}_A)[1]  \buildrel\gamma\otimes \xi\over\lra 
      \bigoplus_{\dim(A)=2} j_{A!} (\ul{\Gc_A}_A)[2] \buildrel\gamma\otimes \xi\over\lra \cdots,
  \]
  the direct sums ranging over all faces of $\Hc$ of given dimension. 
  \end{prop}
 
 \noindent{\sl Proof:} See, e.g., \cite{KS} Eq. (1.12). \qed
 
 \begin{cor}\label{cor:rgamma-c-real}
 If $\Gc\in D^b(V_\RR, \Sc_\RR)$ corresponds to $(\Gc_A,\gamma_{AB})$, then
 \[
 R\Gamma_c(V_\RR, \Gc) \,\,\simeq \,\, \Tot\biggl\{ \bigoplus_{\dim(A)=0} \Gc_A\otimes\OR_A
 \buildrel \gamma\otimes\eps\over\lra \bigoplus_{\dim(A)=1} \Gc_A\otimes\OR_A
  \buildrel \gamma\otimes\eps\over\lra \cdots\biggr\}
 \]
 (the cohomology with compact supports is calculated by the cellular cochain complex). 
  \end{cor}
  
  \noindent {\sl Proof:} This follows because $R\Gamma_c(V_\RR, j_{A!} \ul \k_A) = \OR(A)[-\dim(A)]$
  (cohomology of a cell with compact support). \qed

 \paragraph{B. The complex setup.} 
 Let $V_\CC = \CC^n$ be the complexification of $V_\RR$, and $\Hc_\CC$ the arrangement of hyperplanes in $V_\CC$
 formed by the $H_\CC$, the  complexifications of the hyperplanes $H\in \Hc$. By a {\em flat} of $\Hc_\CC$ 
 we will mean a subspace
 of the form
 $L=\bigcap_{H\in J} H_\CC$ for a subset $J\subset\Hc$ (with  $J= \emptyset$ or $J=\Hc$  allowed). Flats
 form a poset $\Fl(\Hc_\CC)$ ordered by inclusion. Because $\Hc$ is assumed central, $\Fl(\Hc_\CC)$ has
 $0$ as the minimal element and $V_\CC$ as the maximal element. 
 
 For a flat $L$ we denote
 its {\em generic part} by 
  \be\label{eq:gen-part-L}
 L^\circ \,\,=\,\, L  \,\,\setminus\,\,  \bigcup_{H\in\Hc, \,\,  H_\CC \not\supset L} L\cap H_\CC. 
 \ee
The subsets $L^\circ$ form a stratification of $V_\CC$ which we denote by $\Sc_\CC = \Sc_{\CC, \Hc}$.
We view it as a poset, isomorphic to the poset of flats. 
 
  \vskip .2cm
  
  Note that faces can be defined as connected components of $L^\circ_\RR = L^\circ\cap V_\RR$ for strata
   $L^\circ$
  of $\Sc_\CC$. We therefore have the morphism of posets (``complexification")
  \[
  c: \Sc_\RR \lra\Sc_\CC. 
  \]

 We denote by $D^b(V_\CC, \Sc_\CC)$ the full subcategory in the bounded derived category of
 sheaves of $\k$-vector spaces on $V_\CC$ formed by complexes whose cohomology sheaves 
 are constructible with respect to $\Sc_\CC$.  This category has a perfect duality given by passing from
 $\Fc$ to $\Fc^*$, the Verdier dual  complex. 
 Inside it, we have 
  $\Perv(V_\CC, \Sc_\CC)$ the abelian subcategory of perverse sheaves.
     We normalize the conditions of (middle) perversity so that  $\ul\k_{V_\CC}[n]$,
 the constant sheaf put
 in degree $(-n)$,  is perverse.  This normalization agrees with that of \cite{BBD} and differs by shift from
 that of \cite{KS}.  The abelian category  $\Perv(V_\CC, \Sc_\CC)$  is closed under  Verdier duality.  
 
 \paragraph{C. Real data: stalks and hyperbolic stalks.} 
 Let $i_\RR; V_\RR\to V_\CC$ be the embedding. It induces exact functors of triangulated categories
 \[
 i_\RR^*, i_\RR^!: D^b(V_\CC, \Sc_\CC) \lra D^b(V_\RR, \Sc_\RR). 
 \]
 To every complex $\Fc\in D^b(V_\CC, \Sc_\CC)$ and every face $A\in\Sc_\RR$ we can associate therefore 
 two complexes of vector spaces, which we call the {\em stalk} and the {\em hyperbolic stalk}
 of $\Fc$ at $A$:
 \[
 \Fc_A \,=\, (i_\RR^*\Fc)_A \,=\, R\Gamma(A, i_A^* i_\RR^*\Fc), \quad E_A(\Fc) \,=\, 
 (i_\RR^!\Fc)_A \,=\,
 R\Gamma(A, i_A^* i_\RR^!\Fc). 
 \]
 For any pair of faces $A\leq B$ we have the generalization maps (morphisms of complexes) for $i_\RR^*\Fc$ and $i_\RR^!\Fc$:
 \be
 \digamma_{AB}: \Fc_A \lra \Fc_B, \quad \gamma_{AB}: E_A(\Fc)\lra E_B(\Fc). 
 \ee

  By the Duality Theorem, see \cite{KS} Prop. 4.6 or \cite{BFS} Pt. I, Thm. 3.9, we have natural
 isomorphisms
 \be\label{eq:duality}
 E_A(\Fc^*) \,\,\simeq \,\, E_A(\Fc)^*. 
 \ee
 which imply the following.
 
 \begin{prop}\label{prop:stalks and hyperstalks}
 (a) We have a canonical identification
 $E_A(\Fc) \,\simeq\, R\Gamma(A, i_A^! i_\RR^*\Fc)$. 
 
 \vskip .2cm
 
 (b) The hyperbolic stalk $E_A(\Fc)$ is identified with the complex 
 \[
 \Fc_{\geq A}  \,\, :=  \,\, \on{Tot}\, \biggl\{\Fc_A \buildrel\digamma\otimes\eps\over\lra 
 \bigoplus_{B>_1 A} \Fc_B \otimes\OR_{B/A} 
 \buildrel\digamma\otimes\eps\over\lra 
 \bigoplus_{B>_2 A} \Fc_B \otimes\OR_{B/A} \buildrel\digamma\otimes\eps\over\lra \cdots\biggr\}
\]
with the differential $\digamma\otimes \eps$ having matrix elements $\digamma_{BC}\otimes\eps_{BC}$, 
 $B<_1 C$. 
 
 \qed
 
\end{prop}

For a dual statement, expressing ordinary stalks through hyperbolic stalks, see 
Corollary \ref{cor:stalks-and-hyper}.

\vskip .2cm

\noindent {\sl Proof:} Part (a) follows from \eqref{eq:duality} and  the fact that Verdier duality interchanges $i^*$ and $i^!$. 
Part (b) follows by interpreting $i_A^! i_\RR^*\Fc$ as $\ul{R\Gamma}_A(i_\RR^*\Fc)$, 
the complex  of sheaves formed by  (derived) sections with support in $A$. 
The stalk of this complex at any $a\in A$ can be seen as
 \[
 R\Gamma_{\{a\}}(D, i_\RR^*\Fc) \,\,=\,\,R\Gamma_c (D, i_\RR^*\Fc),
\]
  where $D\subset V_\RR$
is a  small transverse open ball  (of complementary dimension) to $A$ centered at $a$. The  situation is similar to that 
of  Corollary \ref{cor:rgamma-c-real}  (with a ball instead of a vector  space) and the same argument gives
the result. 
 \qed

 \vskip .2cm

  It was proved in \cite{KS} Prop. 4.9(a) that for 
 $\Fc\in\Perv(V_\CC, \Sc_\CC)$ the complex $i_\RR^!(\Fc)$ is exact in degrees $\neq 0$, and so the functor
 \be\label{eq:E(F)}
  \Perv(V_\CC, \Sc_\CC) \to \Sh(V_\RR, \Sc_\RR),\quad 
 \Fc\,\mapsto \,  \Ec(\Fc) \, :=\,\ul H^0(i_\RR^!\Fc) \, =\,\ul\HH^0_{V_\RR}(\Fc)
\ee
 is an exact functor of abelian categories. In particular, each $E_A(\Fc)$ reduces to a single vector space. 
   Further, \eqref {eq:duality}
   allows us to define maps of vector spaces
 \[
 \delta_{BA} = 
 \delta_{BA}^\Fc: E_{B}(\Fc)\lra E_A(\Fc), \,\, A\leq B, \quad \delta_{BA}^\Fc := (\gamma_{AB}^{\Fc^*})^*. 
 \]
 which form an anti-representation of $\Sc_\RR$, i.e., a contravariant functor \break
 $(\Sc_\RR,\leq)\to\Vect_\k$. This leads to the following concept.

 \paragraph{D. Hyperbolic sheaves.}
 By a {\em hyperbolic sheaf} on $\Hc$ we will mean a datum
 \[
 \Qc \,=\,   \bigl( E_A, \gamma_{AB}: E_A\to E_{B}, \, \delta_{BA}:
 E_{B}\to E_A, A\leq B \bigr)
\]
 where $E_A, A\in\Sc_\RR$, are finite-dimensional $\k$-vector spaces, $(\gamma_{AB})$ form
 a representation of $\Sc_\RR$, and $(\delta_{BA})$ form an anti-representation so that the following additional
 conditions hold: 
 \begin{itemize}
 \item[(i)] For each $B\leq A$, $\gamma_{BA}\delta_{AB} = \Id_{E_A}$. 
 This allows us to define for  arbitrary $A, B\in  \Sc_\RR$,  
the ``flopping operator''
\[
\phi_{AB} := \gamma_{CB}\delta_{AC}:\ E_A \lra E_B.
\]
Here  $C\in \Sc_\RR$ is any face such that $C\leq A, B$, 
 and the definition  does not depend on the choice of $C$. 
 
 \item[(ii)] Let us call a triple of faces $(A, B, C)$ {\it collinear} if there exist 
points $x\in A, y\in B, z\in C$ lying on the same straight line, 
with $y \in [x, z]$. Then for any such collinear triple we must have
\[
\phi_{AC} = \phi_{BC}\, \phi_{AB}.
\]

\item[(iii)] 
Let $A, B$ be two faces.  Let us say that they are 
{\it neighbors} if they have the same dimension $d$, and  there exists a face 
$C\leq A,  C\leq B$, with $\dim C = d - 1$ 
(a {\it wall} separating $A$ and $B$). Such a wall is unique if it exists. For any such pair of neighbors
we require  that 
$\phi_{AB}$ is an isomorphism. 
 \end{itemize}
 We denote by $\Hyp(\Hc)$ the abelian category formed by hyperbolic sheaves on $\Hc$.
 This category has a perfect duality
 \[
 \Qc=(E_A, \gamma_{AB}, \delta_{BA}) \,\,\mapsto \,\, \Qc^* = (E_A^*, \delta_{BA}^*, \gamma_{AB}^*). 
 \]
  The main result of \cite{KS} can be formulated as follows.
 
 \begin{thm}\label{thm:perv-arr}
 The functor
  \[
\Fc\mapsto  \Qc(\Fc) \,=\, \bigl( E_A(\Fc), \gamma_{AB}: E_A(\Fc)\to E_{B}(\Fc), \, \delta_{BA}:
 E_{B}(\Fc)\to E_A(\Fc), A\leq B \bigr)
 \]
 defines an equivalence $\Perv(V_\CC, \Sc_\CC)\to \Hyp(\Hc)$. This equivalence commutes with duality:
 $\Qc(\Fc^*) \simeq \Qc(\Fc)^*$. 
 \qed
 \end{thm} 
 
The goal of this paper is to describe various features of perverse sheaves explicitly, in terms of
the linear algebra data given by the associated hyperbolic sheaves.

Let us first note the following.

 \begin{prop}\label{f-Rgamma-c}
 If $\Fc\in\Perv(V_\CC, \Sc_\CC)$ corresponds to a hyperbolic sheaf $\Qc(E_A, \gamma_{AB}, \delta_{BA})$,
 then
 \[
 \begin{gathered} 
 R\Gamma_c(V_\CC, \Fc)  \,\,\simeq \,\, \biggl\{ \bigoplus_{\dim(A)=0} E_A\otimes\OR_A
 \buildrel \gamma\otimes\eps\over\lra \bigoplus_{\dim(A)=1} E_A\otimes\OR_ A
  \buildrel \gamma\otimes\eps\over\lra \cdots\biggr\},
  \\
  R\Gamma(V_\CC, \Fc)  \,\,\simeq \,\, \biggl\{ \bigoplus_{\codim(A)=0} E_A\otimes\OR_A
 \buildrel \delta\otimes\eps\over\lra \bigoplus_{\codim(A)=1} E_A\otimes\OR_A
  \buildrel \delta\otimes\eps\over\lra \cdots\biggr\}. 
  \end{gathered}
 \]
 
 \end{prop} 
 
 \noindent {\sl Proof:} The first quasi-isomorphism follows from Corollary \ref{cor:rgamma-c-real}
 and the lemma below. The second quasi-isomorphism follows from the first one by applying
 the Verdier duality. \qed
 
 \begin{lem}
 For any $\Fc\in D^b(V_\CC, \Sc_\CC)$ we have
 \[
 R\Gamma_c(V_\CC, \Fc) \,\,\simeq \,\, R\Gamma_c(V_\RR, i_\RR^!\, \Fc). 
 \]
 \end{lem}
 
 \noindent{\sl Proof of the lemma:}  Let $i_{0, \CC}: \{0\}\to V_\CC$ and $i_{0,\RR}: \{0\}\to V_\RR$
 be the embedings of the origin. Any $\Fc\in D^b(V_\CC, \Sc_\CC)$ is $\RR_+$-{\em conic},
 i..e, each cohomology sheaf of $\Fc$ is locally constant on each orbit of the scaling action
 of $\RR_{>0}$ on $V_\CC$. This implies that 
 \[
 R\Gamma_c(V_\CC, \Fc) \,\simeq \, R\Gamma_{\{0\}}(V_\CC, \Fc) 
 \, =  \, R\Gamma(V_\CC, i_{0,\CC}^!\Fc).
 \]
 Similarly, $i_\RR^!\Fc$ is $\RR_+$-conic on $V_\RR$ and
 \[
 R\Gamma_c(V_\RR, i_\RR^!\Fc)  \,\simeq \, R\Gamma_{\{0\}}(V_\RR, i_\RR^! \Fc) 
  \, =  \, R\Gamma(V_\CC, i_{0,\RR}^! i_\RR^! \Fc),
 \]
 which is the same as the above because $i_\RR i_{0, \RR} = i_{0,\CC}$. \qed
 
 \vskip .2cm
 
 We can now complement Proposition \ref{prop:stalks and hyperstalks} by a
 ``Koszul dual'' statement. 
 
 \begin{cor}\label{cor:stalks-and-hyper}
 For $\Fc\in \Perv(V_\CC, \Sc_\CC)$ the ordinary stalk $\Fc_A$, $A\in\Sc_\RR$ is 
 expressed through hyperbolic stalks as follows:
 \[
 \Fc_A \,\,\simeq \,\, \biggl\{ \bigoplus_{B\geq A\atop
 \codim(B)=0} E_B\otimes\OR_{B/A} \buildrel\delta\otimes\eps\over\lra
 \bigoplus_{B\geq A\atop
 \codim(B)=1} E_B\otimes\OR_{B/A} \buildrel\delta\otimes\eps\over\lra\cdots\biggr\}.
  \]
 \end{cor}
 
 \noindent{\sl Proof:} For $A=0$ this is the second identification of Proposition \ref {f-Rgamma-c},
 since $\Fc_0=R\Gamma(U, \Fc)$ for a small convex open  $U\ni 0$, and this complex
 is independent of $U$, so is the same for $U=V_\CC$. 
 
 For an arbitrary $A$  the statement reduces to the above by considering the quotient arrangement
 $\Hc/L_\RR$ in $V_\RR/L_\RR$,  where $L_\RR$ is the $\RR$-linear span of $A$.
 Faces of  $\Hc/L_\RR$ are in bijection with faces $B$  of $\Hc$ such that $B\geq A$.

 The arrangement  $\Hc/L_\RR$ represents the transversal slice $M$ to $A$; 
 the restriction  $\Fc|_{M_\CC}$ to the complexified
 transversal slice is, by \cite{KS} Prop. 5.3, represented by the hyperbolic sheaf $\Qc^{\geq A}$
 formed by $E_B, B\geq A$, so the calculation of 
 \[
 \Fc_A\, = \,
  R\Gamma(M_\CC, \Fc|_{M_\CC})\,=\, (\Fc|_{M_\CC})_0
  \]  
  reduces to the above case.  
 \vfill\eject

\section {Vanishing cycles in terms of hyperbolic sheaves}

 The standard microlocal approach 
  to study of  perverse sheaves on any stratification is in terms of the local systems
   of vanishing cyclies on the generic parts of conormal bundles to the strata, see
   \cite{MV} \cite{KaSha}.  Our first result provides an explicit description of the fibers
     of these local systems 
   for perverse sheaves from  $\Perv(V_\CC, \Sc_\CC)$.

 \paragraph{A. Background on vanishing cycles.} 
     We recall that for any (polynomial) function $f: V_\CC\to\CC$ and any perverse sheaf $\Fc$
on $V_\CC$ we have a perverse sheaf $\Phi_f(\Fc)$ on $V_\CC$ supported on the hypersurface $\{f=0\}$
and  known
as the {\em perverse sheaf of vanishing cycles}, see \cite{Be}\cite{De}. We will use the 
following real analytic  interpretation of this
perverse sheaf \cite{KaSha}. This interpretation reflects the intuitive meaning of the term
``vanishing cycles". 

\begin{prop}\label{prop:van=rgamma}
We have an isomorphism in the derived category of sheaves on $V_\CC$:
\[
\Phi_f(\Fc) \,\,\simeq \,\, i^*\ul{R\Gamma}_{\{\Re(f)\geq 0\}} (\Fc), 
\]
where $i$ is the closed embedding of the subset $\{f=0\}$ into $\{\Re(f)\geq 0\}$.
\qed
\end{prop}
   
     \vskip .2cm
     
   We will be interested in the case when $f$ is linear. More precisely, 
   let $L^\circ\in\Sc_\CC$ be a stratum, i.e., the generic part of a flat $L$, as in
   \eqref{eq:gen-part-L}. The conormal bundle to $L^\circ$ is
   \[
  T^*_{L^\circ} V_\CC \,=\,  L^\circ \times (V_\CC/L)^* \,\subset \, V_\CC \times V_\CC^* \,=\, T^* V_\CC. 
   \]
   A hyperplane $\Pi\subset V_\CC$ is said to be {\em  transversal to $\Sc_\CC$ at $L$} 
if $L\subset \Pi$, and $L' \in\Fl(\Hc_\CC)$ with $L'\subset \Pi$ implies $L'\subset L$.   
Let us call a {\it polarization at $L$} a linear function 
$f: V_\CC \to \CC$ such that 
$\Pi := \Ker f$ is transversal to $\Sc_\CC$ at $L$.  Polarizations of $L$ form an
open subset  $\Pol(L)\subset (V_\CC/L)^*$, and we define the generic part of the conormal bundle to $L^\circ$ 
as
\[
(T^*_{L^\circ} V_\CC)^\circ \,\,=\,\, L^\circ \times\Pol(L). 
\]

\begin{prop}
Let $\Fc\in\Perv(V_\CC, \Sc_\CC)$. 
If $L\in\Fl(\Hc_\CC)$ and $f\in\Pol(L)$, then $\Phi_f(\Fc)$ is supported on $L$. 
In particular, being perverse, it reduces to a local system in degree $(-\dim(L))$ on $L^\circ$. 
\end{prop}

\noindent{\sl Proof:} 
Let $x\in \{f=0\} \subset V_\CC$ and suppose $x\notin L$. Since  $f\in\Pol(L)$,
  the hyperplane  $\Pi=\{f=0\}$ cannot contain any flats $L'$
  which are not contained in $L$.  So $x$ is not contained in any flat other than $V_\CC$ itself,
  which means that near $x$ the perverse sheaf $\Fc$ is reduced to a local system in degree $(-n)$,
  and so $\Phi_f(\Fc)_x=0$. \qed

\vskip .2cm

 We now describe   the  stalks of the local system $\Phi_f(\Fc)$ at the maximal faces of $L_\RR$.

\paragraph{B.  The complex result.}  

\begin{thm}\label{thm:vanishing}
Let $\Fc\in\Perv(V_\CC, \Sc_\CC)$ and $\Qc = (E_A, \gamma_{AB}, \delta_{BA})
$ be the corresponding hyperbolic sheaf as in
Theorem \ref {thm:perv-arr}.   
  Suppose further that $f\in\Pol(L)$ is real, i.e., takes $V_\RR$ to $\RR$. Let $A$ be a connected
component of $L^\circ_\RR$, so $A$ is
a face of $\Hc$. Consider the complex
\[
E^\bullet_{f, A} \,\,=\,\,\biggl\{ E_A\buildrel\gamma\otimes\eps\over\to
\bigoplus_{B >_1A, \,\,  f|_B \geq 0 } E_B\otimes \OR_{B/A}
\buildrel\gamma\otimes\eps\over \to  \bigoplus_{B >_2 A, \,\,  f|_B \geq 0  } E_B\otimes \OR_{B/A}
\buildrel\gamma\otimes\eps\over \to \cdots \biggr\} 
\]
with the differential $\gamma\otimes\eps$ having matrix elements $\gamma_{BC}\otimes\eps_{BC}$, 
 $B>_1 C$. Then $E^\bullet_{f,A}$
 is exact outside of the leftmost term, and its leftmost cohomology is identified with the  vector space
  $\Phi_f(\Fc)_a[-\dim(L)]$  for any $a\in A$.
 
  \end{thm}
  
  The theorem implies that the shifted space of vanishing cycles  is identified with the subspace
  \[
  E_{f,A} \,\,=\,\, H^0(E^\bullet_{f,A}) \,\,= \bigcap_{B>_1 A, \,\,  f|_B \geq 0 }
  \Ker(\gamma_{AB})\,\,\subset \,\, E_A. 
  \]
 It also implies the following.
  
  \begin{cor}\label{cor:vanishing}
  Consider the complex
 \[
\widecheck E^\bullet_{f, A} \,\,=\,\, \biggl\{ \cdots  \buildrel\delta\otimes\eps\over \to
  \bigoplus_{B >_2 A, \,\,  f|_B \geq 0  } E_B\otimes \OR_{B/A}
   \buildrel\delta\otimes\eps\over \to
     \bigoplus_{B >_1 A, \,\,  f|_B \geq 0  } E_B\otimes \OR_{B/A}
      \buildrel\delta\otimes\eps\over \to
      E_A \biggr\}  
 \]
 with the differential $\delta\otimes\eps$ having matrix elements
  $\delta_{CB}\otimes\eps_{CB}$, 
 $B>_1 C$. Then $E^\bullet_{f,A}$
 is exact outside of the rightmost term, and its rightmost cohomology is identified with the  vector space
  $\Phi_f(\Fc)_a[-\dim(L)]$  for any $a\in A$. In other words, 
   \[
  E_{f,A} \,\,\simeq \,\, \Coker\biggl(  \sum \delta_{BA}: 
   \bigoplus_{B>_1 A, \,\,  f|_B \geq 0 } E_B
    \lra E_A\biggr). 
  \]

  \end{cor} 
  
  \noindent{\sl Proof of the corollary:} The vanishing cycle functor commutes with Verdier
  duality. Therefore the vector spaces $\Phi_f(\Fc)_a[-\dim(L)]$ and 
  $\Phi_f(\Fc^*)_a[-\dim(L)]$ are canonically dual to each other. On the other hand,
  the hyperbolic sheaf corresponding to $\Fc^*$ is, by Theorem \ref{thm:perv-arr},  identified with
  $\Qc^* = (E_A^*, \delta_{BA}^*, \gamma_{AB}^*)$. Our statement follows by
  combining  this
  with Theorem \ref{thm:vanishing} for $\Fc$ and $\Fc^*$. 
  \qed
  
  \begin{rem}
  Theorem \ref{thm:vanishing} and Corollary \ref{cor:vanishing} can be interpreted as follows.
 The same graded space $E^\bullet_{f,A}$ possesses two differentials going in the opposite directions: 
 one induced by the maps $\gamma$, and the other one induced  by the maps  $\delta$. 
 It is natural therefore to form the  ``Laplacian''  $\Delta = \delta\gamma +
\gamma\delta$ out of them.

In the examples we have calculated, $\Delta: E^i_{f,A}\to E^i_{f,A}$ is an isomorphism for $i > 0$. This of course implies the acyclicity statements above.
One may wonder if  this stronger property (Laplacian being an isomorphism for $i>0$) 
holds more generally.  
  \end{rem}

 \paragraph{C. The real analog.} Before proving Theorem \ref{thm:vanishing}, we establish 
 its real counterpart.
 
 \vskip .2cm
 
  Let $\Gc\in D^b(V_\RR, \Sc_\RR)$  and let $(\Gc_A, \gamma_{AB})$
 be the complex of representations of $\Sc_\RR$ corresponding to $\Gc$ by 
 Proposition \ref{prop:constr-gener}.  That is,  $\Gc_A$ is the ordinary stalk of
 $\Gc$ at $A$, and $\gamma_{AB}$ is the generalization map. 
 
 \vskip .2cm
 
 Given a nonzero  $f\in V_\RR^*$, we have the real hyperplane $\Pi =  \{f=0\}\subset V_\RR$. 
 The arrangement $\Hc$ cuts out an arrangement $\Hc\cap\Pi$ in $\Pi$. 
 We denote by $\Sc_{\RR, \Pi}$ the stratification of $\Pi$ into cells of $\Hc\cap\Pi$. 
 We then have the  real version of the  vanishing cycle sheaf.  It is the  complex of sheaves
 \[
i_\Pi^*  \ul{R\Gamma}_{f\geq 0} (\Gc) \,\,\in \,\, D^b(\Pi, \Sc_{\RR, \Pi}). 
 \]
 Here $i_\Pi: \Pi\to V_\RR$ is the embedding. 
 
 \begin{prop}\label{prop:real-vanish} 
  (a)  Let $C'$ be a cell of $\Hc\cap\Pi$ and $C$ be the unique cell of $\Hc$ such that $C'=C\cap\Pi$. The stalk of  $\ul{R\Gamma}_{f\geq 0} (\Gc)$ at $C'$ is quasi-isomorphic to 
  the total
 complex of the  double complex
  \[
  \biggl\{  \Gc_C \buildrel \gamma\otimes\eps\over\lra
 \bigoplus_{D>_1 C,  \,\, f|_D\geq 0} \Gc_D\otimes\OR_{D/C} \buildrel \gamma\otimes\eps \over\lra 
 \bigoplus_{D>_2 C, \,\,  f|_D\geq 0} \Gc_D\otimes\OR_{D/C}  \buildrel \gamma\otimes\eps \over\lra 
 \cdots
 \biggr\} 
 \]
 
\vskip .2cm
 
 (b) Let $C'_1 \leq C'_2$ be an inclusion of cells of $\Hc\cap\Pi$. The generalization map
 \[
 \gamma_{C'_1, C'_2}: \ul{R\Gamma}_{f\geq 0} (\Gc)_{C_1} \lra
 \ul{R\Gamma}_{f\geq 0} (\Gc)_{C_2}
 \]
  is given by the maps $\gamma_{DD'}$ for $\Gc$ which induce  a morphism
 of complexes in (a). 
\end{prop}
 
 \noindent {\sl Proof:} Let $x\in C'$ and $U$ be a small open ball centered at $x$. By definition,
 \[
 \ul{R\Gamma}_{f\geq 0} (\Gc)_{C'} \,\,=\,\, R\Gamma (U, U\cap\{f<0\}; \Gc)
 \]
 The relative cellular  cochain complex representing this, is precisely the complex in (a).
 Part (b) also follows immediately. \qed

 \paragraph{D. Proof of  Theorem \ref{thm:vanishing}.} 
Let $f$ be as in the theorem.   
 Considering $f$ as a complex functional on $V_\CC$, we have
the complex hyperplane $\Pi_\CC = \{f=0\}\subset V_\CC^*$ and
the perverse sheaf $\Phi_f(\Fc)$ on $\Pi_\CC$. By Proposition
\ref{prop:van=rgamma} we can express 
 the hyperbolic stalk of $\Phi_f(\Fc)$
at a cell $C'\in\Sc_{\RR, \Pi}$ as
\[
E_{C'}(\Phi_f(\Fc)) \,=\, (\ul{R\Gamma}_{\Pi_\RR}\,  \ul{R\Gamma}_{\Re(f)\geq 0} (\Fc))_{C'}
\,=\,
( \ul{R\Gamma}_{f\geq 0} \, \ul{R\Gamma}_{V_\RR}(\Fc))_{C'}.
\]
Now, the complex (actually a sheaf) $\Gc = \ul{R\Gamma}_{V_\RR}(\Fc)$ on $V_\RR$ is given by
the stalks $E_B$ and generalization maps $\gamma_{BC}$ from the hyperbolic sheaf $\Qc$.
So applying Proposition \ref{prop:real-vanish}  to  this $\Gc$
and to the cell $C'=A$ as in the formulation of theorem, we get the statement. \qed

\begin{rem} It is worth noticing the following contrast between Proposition  \ref{prop:real-vanish}
and  Theorem \ref{thm:vanishing}. If  $\Gc$
is an arbitrary  sheaf (not a complex)
 on $V_\RR$,  then Proposition  \ref{prop:real-vanish} gives, in general,
a complex with several nontrivial cohomology spaces, because $\ul{R\Gamma}_{f\geq 0}(\Gc)$
need not reduce to a single sheaf.   However, in the case when $\Gc$ has the form
$\Gc=\ul{R\Gamma}_{V_\RR}(\Fc)$
for a perverse sheaf $\Fc\in\Perv(V_\CC, \Sc_\CC)$, this complex  is, by Theorem \ref{thm:vanishing},
 quasi-isomorphic to a single vector space 
in degree $0$.

 A  more immediate instance of such special behavior of  the sheaves 
$\ul{R\Gamma}_{V_\RR}(\Fc)$ can be seen from the property (i) of hyperbolic sheaves in 
\S \ref{sec:real-com}D: 
the condition $\delta_{AB}\gamma_{BA}=\Id$ implies that each $\gamma_{BA}$ is
surjective.
\end{rem}

 \vfill\eject
 
 \section{Specialization and hyperbolic sheaves}\label{sec:spec}
 
 \paragraph{A. Generalities on specialization.} We recall the necessary material from \cite{KaSha} \S 4.1-4.2. 
 Let $X$ be a $C^\infty$-manifold, 
 $M\subset X$ a locally closed submanifold and $T_MX$ the normal bundle to $M$ in $X$.
 Any subset $S\subset X$ gives rise to its {\em normal cone} with center $M$,
  which is a closed subset $C_M S\subset T_MX$ depending only on the closure $\ol S$. 
  We will need the following example. 
  
  \begin{ex}\label{ex:C-subspace}
  Let $X$ be a finite-dimensional $\RR$-vector space and $M\subset X$ be an $\RR$-vector subspace. Then
  $T_MX = M\times (X/M)$. If $S$ is also an $\RR$-vector subspace, then, with respect to the above identification, 
  \[
  C_M(S) \,\,=\,\,  (M\cap S) \times \bigl(S/(M\cap S)\bigr). 
  \]
  \end{ex}
  
  For any complex of sheaves $\Gc\in D^b(\Sh_X)$ we have its {\em specialization} at $M$ which is an 
  $\RR_{>0}$-conic
  complex of sheaves $ \nu_M(\Gc) \in  D^b({T_MX})$. We will  later recall its definition in the case we need.

\paragraph{B. The case of sheaves on arrangements.}
 We will study this construction in two related cases, related to the data of a real arrangement
 $(V_\RR,\Hc)$. 
 
 \vskip .2cm
 
 \noindent \ul{Complex case:} $X=V_\CC$, $M=L_\CC$ a complex flat of $\Hc$
 and $\Gc = \Fc\in\Perv(V_\CC, \Sc_\CC)$ a perverse sheaf smooth with respect to $\Sc_\CC$. 
 
 \vskip .2cm
 
 \noindent \ul{Real case:} $X=V_\RR$, $M=L_\RR$ is a real flat and $\Gc\in D^b(V_\RR, \Sc_\RR)$
 is any complex smooth with respect to the cell decomposition $\Sc_\RR$. 
 
 \vskip .2cm
 
 In each of these cases the normal bundle is itself a vector space:
 \be\label{eq:normal-vector}
 T_{L_\CC} V_\CC \,=\,L_\CC \times (V_\CC/L_\CC), \quad 
  T_{L_\RR} V_\RR \,=\,L_\RR \times (V_\RR/L_\RR). 
 \ee
 The subspace $L_\RR$ carries the {\em induced arrangement} $\Hc\cap L_\RR$ formed by
 the hyperplanes $H\cap L_\RR$ for $H\in\Hc$, $H\not\supset L_\RR$. 
 The quotient space $V_\RR/L_\RR$ carries the {\em quotient arrangement} $\Hc/L_\RR$ formed
 by the hyperplanes $H/L_\RR$ for $H\in\Hc$, $H\supset L_\RR$. We  equip
 $T_{L_\RR} V_\RR$ with the {\em product arrangement} 
 \[
 \begin{gathered}
 \nu_L \Hc \,\,:= \,\,  (\Hc\cap L_\RR) \oplus (\Hc/L_\RR) \,\,=
 \\
 \bigl\{ (H\cap L_\RR) \times V_\RR/L_\RR, \,\, H\not\supset L_\RR\bigr\} \,\,\cup 
 \,\, \bigl\{
 L_\RR \times (H_\RR/L_\RR), \,\, H\supset L_\RR\bigr\}. 
 \end{gathered}
 \]
We have a surjective map $\Hc\to \nu_L(\Hc)$ between (the sets of hyperplanes of) the two
arrangements.  Two hyperplanes $H,  H'$ of $\Hc$ can give the same hyperplane of $\nu_L(\Hc)$,
if $H\cap L_\RR = H'\cap L_\RR$ is the same hyperplane in $L_\RR$. 

\vskip.2cm

We denote by 
\be\label{eq:S-nu-product}
\Sc^\nu_{\RR} = \Sc_{1,\RR} \times\Sc_{2,\RR}, \quad \Sc^\nu_{\CC} = \Sc_{1,\CC} \times\Sc_{2,\CC}
\ee
the stratification 
  of $T_{L_\RR} V_\RR$ by the faces
of $\nu_L(\Hc)$, and  the stratification  of $T_{L_\CC} V_\CC$ by the generic parts of the complex
  flats of $\nu_L(\Hc)$. Here $\Sc_{1,\RR}$ is the stratification of $L_\RR$ by the faces of $\Hc\cap L$,
  while $\Sc_{2,\RR}$ is the stratification of $V_\RR/L_\RR$ by the faces of $\Hc/L$, and
  similarly for $\Sc_{i,\CC}$. 
  
  \begin{prop}\label{prop:spec-smooth}
  (a) If $\Fc\in D^b(V_\CC, \Sc_\CC)$, then 
  $\nu_{L_\CC}\Fc\in D^b(T_{L_\CC}V_\CC, \Sc^\nu_{\CC})$. 
  
  \vskip .2cm
  
  (b)  If $\Gc\in D^b(V_\RR, \Sc_\RR)$, then 
  $\nu_{L_\RR}\Gc\in D^b(T_{L_\RR}V_\RR, \Sc^\nu_{\RR})$. 
  \end{prop}
 
 \noindent{\sl Proof:} We treat only the real case (b), the complex case (a) being
 identical. In the proof we simply write $V$ for the ambient vector space $V_\RR$,
 as well  as $L$ for a real flat and so on.
 We denote by $\SS(\Gc)\subset T^*V$ the microsupport of the complex $\Gc$, and similarly for
 complexes of sheaves on other spaces, see \cite{KaSha} Ch. VI. 
 The statement that $\Gc\in D^b(V, \Sc)$, resp. that $\nu_L(\Gc)\in D^b(T_LV, \Sc^\nu)$, 
  is equivalent to
 \[
 \SS(\Gc) \,\,\subset \,\, \bigcup_{P\in\Fl(\Hc)} T^*_P V, \quad \text{resp.} \quad 
 \SS(\nu_L(\Gc)) \,\,\subset \,\, \bigcup_{Q\in\Fl (\nu_L(\Hc))} T^*_Q (L\times (V/L)). 
 \]
 So we deduce the second inclusion from the first. By Theorem 6.4.1 of \cite{KaSha}, for any
 manifold $X$, a submanifold $M$ and a complex of sheaves $\Gc$ on $X$ we have
 \[
 \SS(\nu_M(\Gc)) \,\subset \, C_{T^*_MX} (\SS(\Gc)) \,\,\subset \,\, T_{T^*_MX} T^*X
 \, \buildrel (!)\over \simeq \, T^*(T_MX). 
 \]
 Here $C_{T^*_MX} (\SS(\Gc))$ is the normal cone to $\SS(\Gc)\subset T^*X$, and the identification
 (!) looks,  in our concrete case, as follows. 
 
 \vskip .2cm
 
 We have $T^*V = V\times V^*$, and $T^*_LV = L\times L^\perp$. Therefore
 \[
 \begin{gathered}
 T_{T^*_LV} T^*V \,\,=\,\, T_{L\times L^\perp} (V\times V^*) \,\,=\,\,
( L\times L^\perp) \times \bigl( (V/L) \times L^*\bigr),
\\ 
T^*(T_LV) \,\,=\,\, T^*\bigl(L\times (V/L)\bigr)\,\,=\,\, \bigl(L\times (V/L)\bigr)\times 
(L^*\times L^\perp),
 \end{gathered}
 \]
 and (!) identifies  factors number 1,2,3,4 of the first product with  factors number 1,4,2,3
 of the second one. 
 
 \vskip .2cm
 
 With this understanding, we need to prove that for any flat $P$ of $\Hc$ the normal cone
 $C_{T^*_LV}(T^*_PV)$ is contained in the union of $T^*_Q\bigl(L\times (V/L))$
 over flats $Q$ of the product arrangement in $L\times (V/L)$. In fact,  it
 is contained in a single $T^*_Q\bigl(L\times (V/L))$, where $Q$ is the product flat
 $(P\cap L)\times \bigl( P/(P\cap L)\bigr )$, as follows from Example \ref{ex:C-subspace}. 
 This finishes the proof of Proposition \ref{prop:spec-smooth}.

 \paragraph{C. Specialization of faces as a continuous map. } 
 Given a face $A$ of $\Hc$, the intersection $\ol A \cap L_\RR$ is the closure of a unique face
 of the arrangement $\Hc\cap L_\RR$ which we denote by $\nu'_L(A)$. Further, the image of
 $A$ in $V_\RR/L_\RR$ is a face of the quotient arrangement $\Hc/L_\RR$ which
 we denote by $\nu''_L(A)$. The pair $\nu_L(A) = (\nu'_L(A), \nu''_L(A))$ is then a face
 of the product arrangement  $\nu_L(\Hc)$ which we call the {\em specialization} of $A$.
 
 \begin{prop}
 The closure of $\nu_L(A)$ is identified with the normal cone $C_{L_\RR}(A)$. Thus $\nu_L(A)$ is the 
 interior (complement of the boundary) of $C_{L_\RR}(A)$. 
 \end{prop} 
 
 \noindent{\sl Proof:} This  is similar to Example  \ref{ex:C-subspace}. \qed
 
\begin{figure}[H]
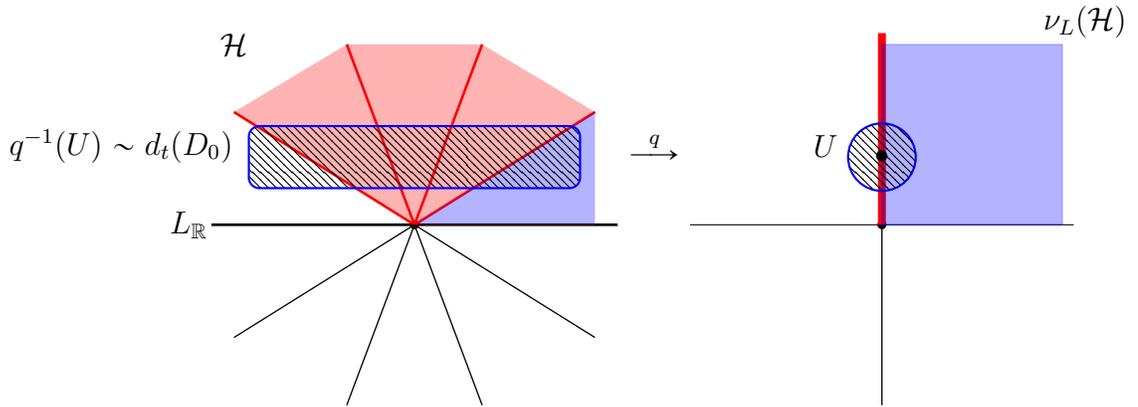
 
  \centering
  
  \btp[scale=.3,
    auto,
    block/.style={
      rectangle,
      draw=blue,
      thick,
      pattern=north west lines, 
      text width=10em,
      align=center,
      rounded corners,
      minimum height=2em
    }]
  
  \node (0) at (0,0){};
  \fill (0) circle (0.2); 
  
  \draw[line width = 1] (-9,0) -- (9,0) {}; 
  \draw [ draw=red, fill=red, opacity = 0.3] (-8,5) -- (0,0)-- (-3,8); 
   \draw [ draw=red, fill=red, opacity = 0.3] (-3,8) -- (0,0)-- (3,8); 
   \draw [ draw=red, fill=red, opacity = 0.3] (8,5) -- (0,0)-- (3,8);   
   \draw[red, line width=1] (-8,5) -- (0,0);   
      \draw[red, line width=1] (-3,8) -- (0,0);   
         \draw[red, line width=1] (3,8) -- (0,0);   
     \draw[red, line width=1] (8,5) -- (0,0);   
   \draw[draw=red, line width=0, fill=blue, opacity=0.3] (8,5) -- (0,0) -- (8,0);   
   \draw [line width = .5] (-8,-5) -- (0,0);      
      \draw [line width = .5] (-3,-8) -- (0,0); 
    \draw [line width = .5] (3,-8) -- (0,0);    
       \draw [line width = .5] (8,-5) -- (0,0);     
  \node at (-10,0) {$L_\RR$};     
  \node at (-8,8){$\Hc$};  
  
    \draw (0,3) node[block](A) {}; 

\node at (-13, 3.5){$q^{-1}(U)\sim d_t(D_0)$}; 
          
  \etp
  \btp
  \node at (0,0){};
 \node at (0,3.3){$\buildrel q\over\lra$};  
  \etp
  \btp[scale=.3]
  
  \node (0) at (0,0){};
  \fill (0) circle (0.2); 
   \draw  (-8.5,0) -- (8.5,0) {}; 
   
   \draw [red, line width=3] (0,0) -- (0,8.5); 
   \draw [line width = 0.5] (0,0) -- (0,-8); 
   \filldraw[color = blue, opacity=0.3] (0,8) -- (0,0) -- (8,0) -- (8,8) -- (0,8); 
   \node at (9,9) {$\nu_L(\Hc)$}; 
  
    \node at (0,3){$\bullet$}; 
    
    \draw [color=blue, thick]  (0,3) circle (1.5cm); 
    
    \fill [pattern=north west lines, 
    ] (0,3) circle (1.5cm); 
    
    \node at (-2.5, 3.5){$U$}; 
            
  \etp

  \caption{Specialization of faces.}\label{fig:spec}

    \end{figure}

 \begin{ex}
 The concept of specialization is illustrated in Fig. \ref{fig:spec}, where $\Hc$ consists of $5$ lines in the
 plane, $L_\RR$ is the horizontal line, and $\Hc/L_\RR$ is the coordinate arrangement
 of two lines in $\RR^2$. The three  open sectors (colored  red) on top, together with
 the open half-lines bounding them, specialize to the upward half-line (also colored   red)
 in $\RR^2$. The open sector (colored blue) with one side being the positive part of $L_\RR$,
 specializes to the first quadrant in $\RR^2$ (also colored blue). 
 \end{ex}
 
    The following is obvious. 
    
     \begin{prop}\label{prop;spec-faces}
  The correspondence $A\mapsto \nu_L(A)$ defines a surjective monotone map
   $\nu_L: \Sc_\RR \to\Sc^\nu_{ \RR}$
  between the posets of faces of $\Hc$ and $\nu_L(\Hc)$ such that
 $\dim \nu_L (A) \leq \dim A$. \qed

  \end{prop}

  We now form the ``geometric realization'' of  the morphism of
  posets $\nu_L$ to construct a continuous map $q: V_\RR\to L_\RR\times(V_\RR /L_\RR)$
  from $V_\RR$ to the normal bundle. 
  That is, 
  choose a point $x_A$ in each face $A$ of $\Hc$. Then we have the
  barycentric subdivision of $V$ into based simplicial convex cones 
  \[
  C(A_1,\cdots , A_p) \,= \, \RR_{>0} \cdot x_{A_1}  + \,\,\cdots \,\,  + \RR_{>0} \cdot x_{A_p}
  \]
   corresponding to all
  increasing chains $A_1<\cdots < A_p$ in $\Sc_\RR$. In particular each $A$ is the union of  the $C(A_1,\cdots, A_p)$ with
  $A_p=A$. 
  Similarly, choose 
 a point $y_B$
  in each face $B$ of $\nu_L(\Hc)$. Then we have the barycentric subdivision of $L\times(V/L)$ into
  similarly defined based  simplicial convex cones $C(B_1,\cdots, B_p)$ for all chains $B_1<\cdots <  B_p$
  in $\Sc^\nu_{\RR}$.    For each chain $A_1<\cdots < A_p$ we 
  define
  \[
  p_{A_1,\cdots, A_p}: \, C(A_1, \cdots, A_p)\lra C(\nu_L(A_1),\cdots, \nu_L(A_p))
  \]
  to be the unique $\RR$-linear map taking $x_{A_i}$ to $y_{\nu_L(A_i)}$. 
  
  \vskip .2cm

  \begin{prop}\label{prop:spec-map}
  $q$ is  a continuous, proper, piecewise linear surjective map. Further, 
  each face $A$ of $\Hc$ is mapped by $q$ to $\nu_L(A)$ in a surjective, piecewise-linear
  way.

  \end{prop}
  
  \noindent {\sl Proof:}  Clear from construction. \qed

  \paragraph{D. The real result.}  In this subsection we deal only with the real situation so
  we write $V$ for $V_\RR$ etc. 
  Let  $\Gc\in D^b(V, \Sc)$
  be a constructible complex.
  
  \begin{thm}\label{thm:spec=direct}
  The specialization  $\nu_L(\Gc)$ is identifed with the
 topological  direct image
 $Rq_*\Gc$ where $q$ is the map from Proposition \ref{prop:spec-map}. 
  \end{thm} 
  
  \noindent{\sl Proof:} We first recall the definition (\cite{KaSha} \S 4.1-2) of $\nu_L(\Gc)$ in terms of the 
  {\em normal deformation} $\wt V_L$ which, in our linear case, reduces to a single chart.
  
  Choose a linear complement $L'$ to $L$ in $V$ so $V=L\oplus L'$.
 Then $L'$ is identified with $V/L$ and $T_LV$ is  also identified with $L\oplus L'$, i.e., with $V$. 
  We  write a general vector
  of $V$ as $v=(l, l')$ with $l\in L$ and $l'\in L'$. Then we define the commutative diagram with Cartesian squares:
  \be\label{eq:blow-diagram}
  \xymatrix{
  T_LV = V\times\{0\} \ar[d]
    \ar[r]^{  s}& \wt V_L: = V\times\RR \ \ar@/^2.0pc/@[][rr]_p
    \ar[d]_{\tau} &\ar[l]_{ \hskip  1cm j}  \Omega \ar[r]^{\wt p}
    \ar[d]^{\wt \tau} & V
  \\
  0 \ar[r] & \RR &\ar[l] \RR_{>0}. &
  }
  \ee
 where  
    \[
  p(l, l', t) = (l, t\cdot l'),\quad 
  \tau  (l, l', t) = t, \quad l\in L, \, l'\in L',\,  t\in\RR. 
 \]
   The space $\Omega$ is defined as $\tau^{-1}(\RR_{>0}) = V\times \RR_{>0}$, and $\wt p$ is the
  restriction of $p$ to $\Omega$. 
  
 After that the specialization is defined by 
  \[
  \nu_L(\Gc)\,\,=\,\, s^* Rj_* \wt p^*(\Gc) \,\, \in
  \,\,  D^b(\Sh_{T_MX}). 
 \]
 Let now $\xi= (l,l')$ be a point of $L\oplus L' = T_LV= \tau^{-1}(0)$. By definition, 
  the stalk of $\nu_L(\Gc)$ at $\xi$ is
  \[
  \nu_L(\Gc)_\xi \,\,=\,\, R\Gamma (D\cap \Omega,  p^*\Gc)
  \]
  where $D\subset V\times\RR$ is a small $(n+1)$-dimensional open ball    around $(\xi,0) = (l, l',0)$.
  Now, $\Omega = V\times\RR_{>0}$. For each $t>0$ consider the slice
  $D_t = D\cap (V\times \{t\})$.  The restriction of $p$ to $D_t$ is the dilation $d_t: (l, l') \mapsto (l, t\cdot l')$
  in the direction of $L'$. 
  
 Since $D$ is a   ball,  the   intersections  $D_t\cap\Omega $ are nonempty for $t$ lying in an open interval of the
 form  $(0,\eps)$ for some $\eps>0$ (the radius of $D$). For such $t$  we have that $D_t\cap\Omega=D_t$
 is the slice over $t$. 
   Since $D$ is a small ball, these
 nonempty slices
  together with the complexes $d_t^*\Gc$ form
  a topologically trivial family over   $(0,\eps)$.  
  This means that we can replace the cohomology of $D\cap\Omega$
  (the union of all slices $D_t, t\in(0,\eps)$) by the cohomology of any single slice, i.e., 
  \[
  \nu_L(\Gc)_\xi \,\,\simeq \,\, R\Gamma( D_t, d_t^*\Gc)
  \]
  for any suffuciently small $t>0$. We can further replace $D_t$ for such $t$ with 
  $0$th slice $D_0 = D\cap (V\times\{0\})$. This slice 
   is  just a small $n$-dimensional open  ball in  $L\oplus L'=V$ around $(l,l')$. This gives
  \[
    \nu_L(\Gc)_\xi \,\,\simeq \,\, R\Gamma( D_0, d_t^*\Gc)\,=\, R\Gamma(d_t(D_0), \Gc), \quad 0<t\ll 1. 
 \]
 When $t\to 0$, the open sets  $d_t(D_0)$ become more and more flattened. 
 We compare them
  with open sets of the form $q^{-1}(U)$ where
 $U$ is a small ball in $T_LV=L\oplus L'$ around $d_t(\xi) = (l, t\cdot l')$. More precisely, we notice that $d_t(D_0)$
 and $q^{-1}(U)$ become homotopy equivalent relatively to the stratification by the faces,
 see Fig. \ref{fig:spec}. This means that we have identifications (the last one expressing the conic nature
 of $Rq_*(\Gc)$: 
 \[
   \nu_L(\Gc)_\xi \,\,\simeq \,\, R\Gamma(q^{-1}(U), \Gc) \,\,=\,\, Rq_*(\Gc)_{d_t(\xi)} \,\,\simeq\,\,
   Rq_*(\Gc)_\xi. 
 \]
 This identifies the stalks. The same considerations show that the generalization maps between the stalks match as well. 
 The theorem is proved. \qed
 
 \vskip .2cm
 
 Assume now that $\Gc$ is given by
  a  complex of representations $G=(\Gc_A, \gamma_{AA'})$ of $\Sc_\RR$.  So the complexes $\Gc_A$ are the stalks of $\Gc$ and
 the  $\gamma_{AA'}$ are the generalization maps. For any face $B\in \Sc^\nu_{\RR}$   of $\nu_L(\Hc)$ define a complex 
      \be
  \Gc_{L,B} \,\, =  \,\,\on{Tot} \,\biggl\{ \bigoplus_{\nu_L(A)=B \atop
 \dim(A)=\dim(B)} \Gc_A\otimes \OR_{A/B}  \buildrel \gamma\otimes\eps \over\lra 
 \bigoplus_{ \nu_L(A)=B\atop
 \dim(A)=\dim(B)+1} \Gc_A \otimes\OR_{A/B}
 \buildrel \gamma\otimes\eps \over\lra 
  \cdots
 \biggr\}. 
 \ee
  Let $B <_d B'$ be two faces of $\nu_L(\Hc)$. 
 We define a morphism of complexes
\[
\gamma^L_{B, B'}: \Gc_{L,B}\to \Gc_{L, B'}
\] 
as follows. Let $A \in\Sc_\RR$ be such that $\nu_L(A)=B$ and $\dim(A)=\dim(B)+p$, 
so that $\Gc_A\otimes \OR_{A/B}$ is a summand in the $p$th term of $\Gc_{L,B}$. 
Similarly  let $A '\in\Sc_\RR$ be such that $\nu_L(A')=B'$ and $\dim(A')=\dim(B')+p$, 
so that $\Gc_A'\otimes \OR_{A'/B'}$ is a summand in the $p$th term of $\Gc_{L,B'}$. 
If $A\leq A'$, 
 then $A< _d A'$ and the identification
  of the quotient spaces
   \[
  \Lin_\RR(A')/\Lin_\RR(A) \buildrel \simeq\over\lra \Lin_\RR (B') /\Lin_\RR(B)
  \]
  gives, passing to the determinants and transposing,  an isomorphism
 \[
 \sigma^*_{AA'}: \OR_{A/B} \lra \OR_{A'/B'}. 
 \]
 We define the matrix element 
 \[
 (\gamma_{B, B'}^L)_A^{A'}: \Gc_A \otimes \OR_{A/B} \lra \Gc_{A'}\otimes\OR_{A'/B}
 \]
 to be equal to $\gamma_{AA'}\otimes\sigma^*_{AA'}$ if $A<A'$ and to $0$ otherwise.

 \begin{cor}\label{cor:real-spec-gamma}
   Each $\gamma^L_{BB'}$ is indeed a morphism of complexes, and the data
 $(\Gc_{L,B}, \gamma^L_{BB'})$ is a complex of representations of $\Sc_{\nu,\RR}$,
 the poset of faces of the arrangement $\nu_L(\Hc)$. 
  This complex of representations  describes the constructible complex $\nu_L(\Gc)$.
\end{cor}

\noindent{\sl Proof:} Choose any point $b\in B$. Since $q$ is a proper map, the stalk of $Rq_*(\Gc)$ at $b$ is identified with
$R\Gamma(q^{-1}(b), \Gc)$. Now
$\Gc_{L,B}$ is nothing but  the cellular cochain complex calculating    $R\Gamma(q^{-1}(b), \Gc)$.
We similarly identify the generalization maps. \qed

\begin{rem}
At the formal algebraic level, 
the property that $\gamma^L_{BB'}$ is indeed a morphism of complexes,
simply reflects the fact that the differential in $R\Gamma(V, \Gc)$, the cellular cochain
complex, satisfies $d^2=0$. More precisely, we have an identification
(isomorphism, not just a quasi-isomorpism) of cellular cochain complexes
\[
R\Gamma(V, \Gc) \,\,\simeq \,\, R\Gamma(L\times (V/L), Rq_*(\Gc)) \,\,
\ \simeq \,\, R\Gamma(L\times (V/L), \nu_L(\Gc)).
\]
The RHS of this identification represents the same complex in a ``block'' form,
with blocks (stalks of $\nu_L(\Gc)$) parametrized by faces $B$ of $\nu_L(\Hc)$.
The fact that the maps $\gamma^L_{BB'}$ between the blocks are morphisms of
complexes is implied by the fact that the total differential squares to $0$. 

\end{rem}

   \paragraph{E. Bispecialization.} 
   We first  consider the general situation  studied in \cite{schapira-takeuchi}
   \cite{takeuchi}. 
   Let $N\subset M\subset X$ be a flag of $C^\infty$ submanifolds in a $C^\infty$ manifold $X$. 
   In the normal bundle $T_N X$ we have the submanifold  (subbundle) $T_N M$. 
   In the normal bundle $T_MX$ we have the submanifold $N$, emdedded into $M$
   (the zero section of $T_MX$). It turns out that the normal bundles of these new
   submanifolds are identified.  
   
   \begin{prop}\label{prop:two-tangents}
     We have identifications\footnote{The notation $\oplus$ here and below means direct sum
       of vector bundles, i.e.\ fiber product over $N$.}
   \[
   T_{T_N M } (T_NX) \,\buildrel (1) \over \simeq \, T_N M \oplus (T_M X)|_N  \,
   \buildrel (2) \over \simeq \, 
   T_N (T_MX). 
   \]
   \end{prop}
   
   \noindent {\sl Proof:} The statement is a part of  Prop. 2.1 of \cite{takeuchi}. For convenience of the reader
   we give a sketch of the proof. The identification 
    (1)  is a particular case of the well known fact which generalizes,
   to vector bundles, 
    the identification \eqref{eq:normal-vector} for vector spaces: If $L\subset V$ is a 
    $C^\infty$ vector
    subbundle in a $C^\infty$ vector bundle over a $C^\infty$-manifold $B$, then
    $T_LV \simeq L\oplus (V/L)$. To see (2),  we recognize, inside $T_N(T_MX)$
    two subbundles: first, $T_NM$ (the normal bundle to $N$ inside the zero section of $T_MX$),
    and, second $(T_MX)|_N$ (the restriction to $N$ of the normal bundle). Inspection in
    local coordinates shows that these two subbundles form a direct sum decomposition.\qed
    
    \vskip .2cm
    
    In this context Schapira and Takeuchi  \cite{schapira-takeuchi}
   \cite{takeuchi} defined a functor
   \[
   \nu_{NM}: D^b (X) \lra D^b ({T_N\oplus (T_XM)|_N})
   \]
   called {\em bispecialization}. It  is defined, similarly to the usual specialization,  through the
   {\em binormal deformation} $\wt X_{NM}$, recalled below. 
      On the other hand, we can iterate the specialization functors, getting a
     diagram
    of  functors between derived categories of sheaves on the manifolds in question:
    \be\label{eq:diag-iter-spec}
      \xymatrix{
   D^b (X)
   \ar[d]_{\nu_M}  \ar[r]^{\nu_N} 
   \ar[dr]^{\nu_{NM}}
   & D^b ({T_NX})
   \ar[d]^{\nu_{T_NM} }
   \\
   D^b({T_MX}) \ar[r]_{\hskip -1cm \nu_N} &D^b ({T_NM \oplus (T_MX)|_N}).  
      }  
    \ee
   This diagram is not (2-)commutative, i.e., the two composite functors  (iterated specializations)
 are not isomorphic.
       
    \begin{ex}
  Let $X=\RR^2$ with coordinates $x,y$, let $M$ be the line $y=0$ and
    $N$ be the origin $(0,0)$.   Let $P\subset X$ be the parabola $y=x^2$
    and $\Gc=\ul\k_P$ be the constant sheaf on $P$. We identify all three manifolds
     $T_NX$, $T_MX$ and
    $T_NM \oplus (T_MX)|_N$  back with $\RR^2$ with the same coordinates. 
     Then $\nu_N(\Gc)$ is the constant sheaf on the horizontal  line $y=0$ (the tangent line to $P$),
     and $\nu_{T_NM}(\nu_N(\Gc))$ is again the constant sheaf on the line $y=0$. 
     On the other hand, $\nu_M(\Gc)$ is supported on the vertical half-line $x=0, y\geq 0$
     (since $P$ is contained in the upper half plane $y\geq 0$ and does not meet $M$
     except for $x=0$). 
     So $\nu_N(\nu_M(\Gc))$ will be
     again supported on this half-line. 
     \end{ex} 
     
     Nevertheless, in the linear case   all three possible functors are identified.

     \begin{thm}\label{thm:bispec}
     Let $X=V$ be an $\RR$-vector space and $N\subset M\subset V$ be a flag of $\RR$-linear
     subspaces. Let $\Hc$ be an arrangement of hyperplanes in $V$ and $\Sc_\RR$
     the corresponding stratification by faces. Then for
     $\Gc\in D^b (V,\Sc_\RR)$
     we have canonical quasi-isomorphisms
    \[
    \nu_N(\nu_M(\Gc))\,\,\simeq \,\, \nu_{T_NM}(\nu_N(\Gc))
    \,\,\simeq \,\,
     \nu_{NM}(\Gc).
     \]
       In other words, the diagram \eqref{eq:diag-iter-spec}
     becomes 2-commutative if the top left corner is replaced by $D^b (V,\Sc_\RR)$. 
     \end{thm}
     
     \noindent{\sl Proof:} Enlarging $\Hc$ if necessary, we can assume that $N$ and $M$
     are flats of $\Hc$. The space $T_NM\oplus (T_MX)|_N$ is identified with
     vector space $V''= N\oplus (M/N) \oplus (V/M)$ which carries the {\em triple product
     arrangement} 
     \[
     \nu_{NM}(\Hc)\,: =\, (\Hc\cap N) \oplus( (\Hc\cap M)/N) \oplus (\Hc/M). 
     \] 
     Denote by $\Sc^{\nu_{NM}}_\RR$ the stratification given by the faces of this arrangement.
     Also denote $\Sc^{\nu_N}_\RR$ and $\Sc^{\nu_M}_\RR$ the stratifications given by the faces
     of the arrangements $\nu_N(\Hc)$ and $\nu_M(\Hc)$. Now notice that specialization of faces gives
     a {\em commutative} diagram of morphisms of posets which we then use to construct a commutative
     diagram of proper piecewise linear maps: 
     \[
     \xymatrix{
     \Sc_\RR\ar[d]_{\nu_M}  \ar[r]^{\nu_N} & \Sc_\RR^{\nu_N}\ar[d] ^{\nu_{T_NM}}
     \\
     \Sc_\RR^{\nu_M} \ar[r]_{\nu'_N} & \Sc_\RR^{\nu_{NM}}
     }
     \quad\quad
     \xymatrix{
     V \ar[r]^{q_N}
     \ar[d]_{q_M}  & T_NV \ar[d]^{q_{T_NM} }
     \\
     T_MV \ar[r] ^{\hskip -1cm q'_N}& T_NM \oplus (T_MV)|_N. 
     }
   \] 
   The direct images in this second diagram correspond, by Theorem \ref{thm:spec=direct},
   to the specialization functors on the outer edges of the diagram \eqref{eq:diag-iter-spec}.
   This shows that the outer rim of \eqref{eq:diag-iter-spec} is 2-commutative. 
   
   We now show that the composite functor given by the outer rim of \eqref{eq:diag-iter-spec}, 
   is isomorphic to $\nu_{NM}$. (This will also  give another  proof of the commutativity of the outer rim.)
   For this we recall the explicit form of the binormal deformation diagram, see
   \cite{takeuchi} Eq. (2.20). We choose a complement $L'$ to $N$ in $M$ and a complement $L''$ to $M$ in $V$,
   thus identifying $V$, as well as $T_NM \oplus (T_MV)|_N$, with $N\oplus L'\oplus L''$. 
   So we write  elements of either of this spaces as $(n, l', l'')$. 
   Then the ``bi''-analog of the diagram 
   \eqref {eq:blow-diagram} has the form 
   
  \[
  \xymatrix{
  T_NM \oplus (T_MV)|_N  =  V\times \{(0,0)\}\ar[d]
    \ar[r]^{\hskip 1.5cm  s}& \wt V_{NM}  = V\times \RR^2 \ar@/^2.0pc/@[][rr]_p
    \ar[d]_{\tau} &\ar[l]_{\hskip 1cm j} \Omega \ar[r]^{\wt p}
    \ar[d]^{\wt \tau} & V
  \\
  0 \ar[r] & \RR^2 &\ar[l] \RR^2_{>0}, &
  }
  \]
  with
  \[
  p((n, l', l''), (t',t'')) \,\,=\,\, (n, \, t'l',  \, t' t'' l''), \quad \tau ((n, l', l''), (t',t'')) = (t', t''),
  \]
  so the restriction of $p$ to $V\times \{(t',t'')\} $ is the map
  \[
  p_{(t',t'')}: (n, l', l'') \mapsto (n,\,  t'l', \,  t' t'' l''). 
  \]
The bispecialization is defined as 
$\nu_{NM}(\Gc) = s^* Rj_* \wt p^*\Gc$ with respect to this diagram, so its stalk at $(n,l',l'')$
is $R\Gamma(D\cap\Omega, p^*\Gc)$ where $D$ is a small open $(n+2)$-dimensional ball around $((n,l',l''),(0,0))$. 
We slice $D$ into $n$-dimensional balls  $D_{(t',t'')} =  D\cap \tau^{-1}(t', t'')$. 

\begin{lem} For sufficiently small $\eps>0$, 
the slices $D_{(t',t'')}$ together with
 the restrictions $p^*\Gc|_{D_{(t',t'')}} = p_{(t', t'')}^*\Gc$, form a topologically trivial family
 over the product of open intervals $(0,\eps) \times(0,\eps)$.  
\end{lem}

\begin{figure}[H]
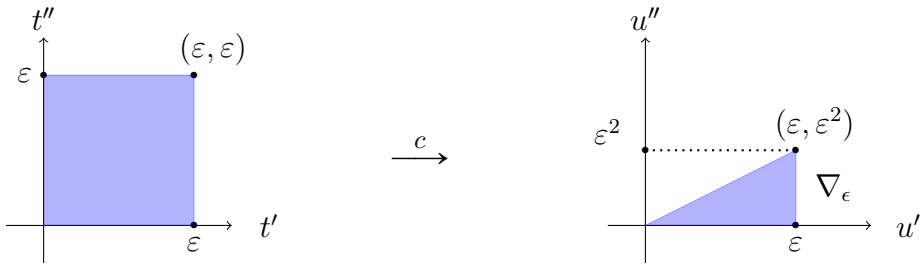
 
  \centering
  
  \btp[scale=.5]
  
  \node at (0,0){\large$\buildrel c\over\lra$};
  
 \draw[->]  (5, -2) -- (12, -2); 
  \draw [->] (6, -3) -- (6,3); 
  
    \node at (10, -2){\tiny$\bullet$}; 
  \node at (6,0){\tiny$\bullet$}; 
 \node at (10,0){\tiny$\bullet$}; 
 
   \filldraw[color = blue, opacity=0.3] (6, -2) -- (10, -2) -- (10,0) -- (6, -2); 
   \draw[dotted, line width =0.8] (6,0) -- (10,0); 
   
   \node at (10, -2.5){$\eps$}; 
   \node at (10.5, 0.7){$(\eps, \eps^2)$}; 
  \node at (5, 0.5){$\eps^2$}; 
  
   \draw[->]  (-11, -2) -- (-5, -2); 
 \draw[->]  (-10, -3) -- (-10, 3);   
  
     \node at (-6, -2){\tiny$\bullet$};  
    \node at (-6,2){\tiny$\bullet$};   
     \node at (-10, 2){\tiny$\bullet$}; 
 
  \filldraw[color = blue, opacity=0.3] (-10, -2) -- (-6, -2) -- (-6,2) -- (-10, 2) -- (-10, -2); 
  
 \node at (-6, -2.5){$\eps$}; 
  \node at (-10.5, 2){$\eps$}; 
\node at (-5.5, 2.7){$(\eps, \eps)$};

\node at (-4, -2){$t'$}; 
\node at (-10, 3.5){$t''$}; 

\node at (13, -2){$u'$}; 
\node at (6, 3.5){$u''$};

\node at (11,-1){$\nabla_\epsilon$}; 
              
  \etp

  \caption{The wedge $\nabla_\epsilon$.}\label{fig:wedge}

    \end{figure}

\noindent{\sl Proof of the lemma:} For $u', u'' >0$ let
\[
d_{(u', u'')}: V\to V, \quad (n, l', l'') \mapsto (n, \, u'l',\,  u'' l'')
\]
be the bi-dilation in the last two variables.
 Then $p_{(t', t'')} = d_{c(t', t'')}$, where $c: \RR^2\to \RR^2$ is the  map
\[
(t', t'') \mapsto (u', u'') = (t',\, t't''). 
\]
Now, $c$ maps the open square  $(0,\eps)^2$ homeomorphically onto the open triangular wedge $\nabla_\eps$
of slope $\eps$, see Fig. \ref{fig:wedge}.

    For small $t', t''>0$ we can identify the slices $D_{(t', t'')}$ with $D_{(0,0)}$ (alternatively, we
    could have taken $D$ to be the product of balls in $V$ and in $\RR^2$ so that the slices 
    would not change at all). 
    
    We recall that $\Gc$ is smooth with respect
    to a hyperplane arrangement $\Hc$ (so the slopes of the hyperplanes are fixed). 
    On the other hand,  the slope of the wedge $\nabla_\epsilon$ is shrinking as $\eps\to 0$.
   Therefore,   for sufficiently small
    $\eps$ we will have that for all $(u', u'')\in\nabla_\epsilon$ the topological structure of   $d_{(u', u'')}^*\Gc$ on $D_{(0,0)}$  will stabilize.  This proves the lemma. \qed
 
 \vskip .2cm
 
 The lemma implies that
   the stalk of $\nu_{NM}(\Gc)$  at $(n,l',l'')$ can be written as 
 \[
 R\Gamma (D_{(0,0)}, p_{(t',t'')}^*\Gc)\,=\, R\Gamma(p_{(t', t'')}(D_{(0,0)}), \Gc)
 \]
  for any sufficiently small
 positive $t', t''$.
 
 It remains to  similarly analyze the two outer composite functors (iterated specializations) in 
 \eqref{eq:diag-iter-spec}  and to  find that they correspond to the choice of $0< t' \ll t''\ll 1 $, resp. $0< t'' \ll t' \ll 1$. Because
 of the topological triviality of the family over all $(t', t'')\in (0,\eps)\times (0,\eps)$, all three  results are the same. 
 \qed

 \paragraph{ F. The complex result.} 
 We now consider the complex situation: that of a perverse sheaf 
 $\Fc\in\Perv(V_\CC, \Sc_\CC)$ and the corresponding 
  hyperbolic sheaf $\Qc=(E_A, \gamma_{AA'},\delta_{A'A})$. 
   Let $\Qc^\nu = (E^\nu_B, \gamma^\nu_{BB'}, \delta^\nu_{B'B})$ be the hyperbolic sheaf corresponding
 to $\nu_{L_\CC}(\Fc)\in \Perv(T_{L_\CC} V_\CC, \Sc_{\nu , \CC})$. Here $B, B'$ are faces of the product
 arrangement $\nu_L(\Hc)$.

 \begin{thm}\label{thm:spec-complex}

 (a) The hyperbolic stalk $E^\nu_B$ is identified as 
 \[ 
 E_B^\nu  \,\,\simeq \,\, \,\biggl\{ \bigoplus_{\nu_L(A)=B\atop
 \dim(A)=\dim(B)} E_A\otimes \OR_{A/B}  \buildrel \gamma\otimes\eps \over\lra 
 \bigoplus_{\nu_L(A)=B\atop
 \dim(A)=\dim(B)+1} E_A \otimes\OR_{A/B}
 \buildrel \gamma\otimes\eps \over\lra 
  \cdots
 \biggr\}
 \]
 That is,  the complex in the RHS is exact everywhere except the leftmost term, 
 where the kernel is
 identified with $E_B^\nu$. 
 
 \vskip .2cm
 
 (a') We also have an identification
 
 \[ 
 E_B^\nu  \,\,\simeq \,\, \,\biggl\{    \cdots
  \buildrel \delta\otimes\eps \over\lra 
  \bigoplus_{\nu_L(A)=B \atop
 \dim(A)=\dim(B)+1} E_A \otimes\OR_{A/B}
   \buildrel \delta\otimes\eps \over\lra 
     \bigoplus_{\nu_L(A)=B \atop
 \dim(A)=\dim(B)} E_A \otimes\OR_{A/B}
 \biggr\}. 
 \]
  That is,  the complex in the RHS is exact everywhere except the rightmost term, 
 where the cokernel is
 identified with $E_B^\nu$.
 
  \vskip .2cm
  
   (b) The maps $\gamma^\nu_{BB'}$   are induced by the maps $\gamma_{AA'}$  which induce morphisms of complexes in (a), similarly to Corollary \ref{cor:real-spec-gamma}.
   
     \vskip .2cm

 (b') The maps   $\delta^\nu_{B'B}$  are induced by the map  $\delta_{A'A}$ which induce morphisms of complexes in (a').

 \end{thm}
 
 \noindent {\sl Proof:} We first prove parts (a) and (b).  Let $i_\RR: V_\RR\to V_\CC$ and $i_{\RR, \nu}: T_{L_\RR} V_\RR \to T_{L_\CC} V_\CC$ be the
 embeddings of the real parts.  Put
 \[
 \Gc = i_\RR^!\Fc, \quad \Gc_\nu = i_{\RR, \nu}^!\,  \nu_{L_\CC} (\Fc). 
 \]
 These are ordinary sheaves (not just complexes) on  $V_\RR$ and $T_{L_\RR} V_\RR$,  smooth with respect to
 $\Sc_\RR$ and $\Sc^\nu_{\RR}$ respectively. Their stalks are given by the 
  $E_A$ and $E^\nu_B$  and  their generalization maps are given by the $\gamma_{AA'}$ and $\gamma^\nu_{BB'}$
  respectively. Note that we have a canonical morphism
  \[
  \nu_{L_\RR}(\Gc) \,=\, \nu_{L_\RR} (i_\RR^! \Fc) \buildrel\beta\over\lra i_{\RR,\nu}^! \nu_{L_\CC}(\Fc) \,=\,\Gc_\nu,
  \]
 see \cite{KaSha} Prop. 4.2.5. So our statements will follow from   
 Corollary \ref{cor:real-spec-gamma} if we establish the
 following.
 
 \begin{prop}\label{prop:beta-c}
 For any $\Fc\in D^b(V_\CC, \Sc_\CC)$, the morphism $\beta:  \nu_{L_\RR} (i_\RR^! \Fc)\to i_{\RR,\nu}^! \nu_{L_\CC}(\Fc)$
  is a quasi-isomorphism. 
 \end{prop}
   
   \noindent {\sl Proof of  Proposition \ref{prop:beta-c}:} Since $\nu_{L_\RR}$ and
   $\nu_{L_\CC}$ commute with Verdier duality, it is enough to show that for any 
  $\Fc\in D^b(V_\CC, \Sc_\CC)$, the dual morphism 
  $\alpha: i_{\RR,\nu}^* \nu_{L_\CC}(\Fc) \to \nu_{L_\RR}(i_\RR^*\Fc)$ is a quasi-isomorphism.
   Such a morphism is defined for any $\Fc\in D^b(\Sh_{V_\CC})$
   whatsoever, see \cite{KaSha} Prop. 4.2.5. 
   So we show that it is a quasi-isomorphism for a more general class of complexes.
   Namely, $V_\CC$ has the {\em product stratification} $\Sc_\RR\times \Sc_\RR$ formed by the cells of the form
   $A'+iA''\subset V_\CC = V_\RR + i V_\RR$, where $A'$ and $A''$ are arbitrary faces of $\Hc$ and $i=\sqrt{-1}$. 
   This stratification refines $\Sc_\CC$, so $D^b(V_\CC, \Sc_\CC) \subset D^b(V_\CC, \Sc_\RR\times\Sc_\RR)$. 
   Therefore it suffices to prove:
   
   \begin{lem}\label{lem:alpha-rr} 
   For any $\Fc\in D^b(V_\CC, \Sc_\RR\times\Sc_\RR)$, the morphism 
    $\alpha: i_{\RR,\nu}^* \nu_{L_\CC}(\Fc) \to \nu_{L_\RR}(i_\RR^*\Fc)$
   is a quasi-isomorphism. 
   \end{lem} 
   
    \noindent {\sl Proof of  Lemma \ref{lem:alpha-rr}:} The stratification on $V_\RR$
    induced by $i_\RR$ from $\Sc_\RR\times\Sc_\RR$, is $\Sc_\RR$. This means that
    the specializations maps  of the posets of faces are compatible, and therefore we have a
    commutative diagram
    \[
    \xymatrix{
    V_\RR \ar[r]^{\hskip -1cm q_\RR}
    \ar[d]_{i_\RR}
    & L_\RR\times (V_\RR/L_\RR)\ar[d]^{i_{\RR,\nu}}
    \\
      V_\CC \ar[r]_{\hskip -1cm q_\CC}& L_\CC\times (V_\CC/L_\CC),
      }
    \]
    where $q_\RR$ and $q_\CC$ are the proper maps constructed in
    Proposition \ref{prop:spec-map}. So our statement follows from
    Theorem \ref{thm:spec=direct} by proper base change.\qed
    
    This finishes the proof of Proposition \ref{prop:beta-c} and of parts (a) and (b)  of
    Theorem \ref{thm:spec-complex}. 
    
    \vskip .2cm
    
    Now, parts (a') and (b')  of
    Theorem \ref{thm:spec-complex} follow from (a) and (b) because $\nu_{L_\CC}$ commutes
    with Verdier duality whose effect on hyperbolic sheaves exchanges $\gamma$
    and $\delta$, see Theorem \ref{thm:perv-arr}. 
     Theorem \ref{thm:spec-complex}  is proved. 
    \vfill\eject

   \section{Fourier transform and hyperbolic sheaves}\label{sec:fourier-hyp}
   
   \paragraph{A. Generalities on the Fourier-Sato transform.}  Let $W$ be a finite-dimensional $\RR$-vector
   space and  $W^*$ the dual space. We denote by $D^b_{\con}(E)\subset D^b(E)$ the
   full subcategory formed by complexes $\Gc$ which are {\em conic}, i.e., such that each
   sheaf $\ul H^j(\Gc)$ is locally constant on any orbit of the scaling action of $\RR_{>0}$ on $W$. 
   
   \vskip .2cm
   
   Set
   \[
   P \,\,=\,\,\bigl\{ (x, f) \in W\times W^* \,\bigl| \, f(x)\geq 0\bigr\}\,\, \buildrel i_P\over\hookrightarrow \,\,
   W\times W^*
   \]
and denote 
 by $p_1, p_2$ the projections of $P$ to $W$ and
   $W^*$ respectively.    
   The {\em Fourier-Sato transform} is an equivalence of categories
   \[
   \FS: D^b_\con(W)\lra D^b_\con(W^*), \quad \FS(\Gc) \,=\, Rp_{2!}(p_1^*\Gc),  
   \]
  see \cite{KaSha} Def. 3.7.8.   The base change theorem implies at once the following.
  
  \begin{prop}\label{prop:FS-stalk-K}
  Let $f\in W^*$.  The stalk of $\FS(\Gc)$ at $f$ is found as
  \[
  \FS(\Gc)_f \,\,\simeq \,\, R\Gamma_c(P_f,\Gc), 
  \]
where $P_f= p_2^{-1}(f) = \{x\in W| f(x)\geq 0\}$. (Thus $P_f$ is a closed half-space
  for $f\neq 0$ and $P_f=W$ for $f=0$.)   
  
   \qed
  
  \end{prop}

   \paragraph{B. The dual arrangement.} 
   We specialize the above to the two situations related to 
 an arrangement of hyperplanes $\Hc$ in $V_\RR$. We denote $n=\dim_\RR V_\RR$. 
 
 \begin{itemize}
 \item[(1)] $W=V_\RR$ and $\Gc\in D^b(V_\RR, \Sc_\RR)$. In this case we would like to
  find the stalks of $\FS(\Gc)$.
 
 \item[(2)] $W=V_\CC$ and $\Gc \in \Perv(V_\CC, \Sc_\CC)$. 
 We identify $W^*=\Hom_\CC(V,\CC)$ with the real dual $\Hom_\RR(V_\CC, \RR)$
 by means of the form 
 \[
 (x,f)\mapsto \Re (f(x)), \quad x\in V_\CC, f\in V_\CC^*.
 \]
 
 In this case it is known, see \cite{KaSha} Ch. X, 
 that $\FS(\Gc)[-n]$ is a perverse sheaf on $V_\CC^*$ with respect to some stratification.
 We would like to relate this stratification to an arrangement of hyperplanes and to find
 the hyperbolic stalks of $\FS(\Gc)$. 
 \end{itemize}
 
 This leads to the following definition.
 
 \begin{Defi}
 The {\em dual arrangement} $\Hc^\vee$
   of hyperplanes in $V_\RR^*$ consists of orthogonals $l^\perp$
   where $l$ is a 1-dimensional flat of $\Hc$. We denote by $\Sc^\vee_\RR$ the stratification
   of $V_\RR^*$ into faces of $\Hc^\vee$ and by $\Sc^\vee_\CC$ the stratification of $V_\CC^*$ into
   generic parts of the complex flats of $\Hc^\vee$. 
   
   \end{Defi} 
   
   \begin{prop}\label{prop:double-dual}
  We have an inclusion $\Hc\subset \Hc^{\vee\vee}$ (as sets of hyperplanes in $V_\RR$).

   \end{prop}
   
   \noindent {\sl Proof:}  1-dimensional flats of $\Hc^\vee$ are the
   orthogonals $M^\perp$, where $M$ runs over  hyperplanes in  $V_\RR$
   which are sums of 1-dimensional flats of $\Hc$. Such $M$ are therefore, precisely
   the hyperplanes of $\Hc^{\vee\vee}$. Now the statement  means that each
   hyperplane $H\in\Hc$ can be obtained as a sum of 1-dimensional flats of $\Hc$.
   This is indeed the case, since we have assumed from the outset that $\Hc$
   is central, i.e., the intersection of all $H\in\Hc$ is $0$.  \qed

    \begin{exas}
    (a) Call an arrangement $\Hc$ {\em reflexive}, if $\Hc^{\vee\vee}=\Hc$. A sufficient condition for this is
    that the set of  flats of $\Hc$ is closed not only under intersections but also under sums, i.e., it forms
    a  lattice. This follows from the proof of Proposition \ref{prop:double-dual}. 
    Examples of reflexive  arrangements include any arrangement with $\dim(V_\RR)\leq 2$, as well
    as any direct sum of such arrangements. 
    
    
    \vskip .2cm
    
      (b) 
    In general, forming the  union of the arrangements
   \[
   \Hc \,\subset \, \Hc^{\vee\vee} \,\subset \, \Hc^{\vee\vee\vee\vee} \,\subset \,\cdots
   \]
   amounts to closing $\Hc$ under the operations of sum and intersection, i.e.,
   to forming the  lattice  of subspaces
   generated by $\Hc$ and taking all $(n-1)$-dimensional
   elements of it.  
    Such a lattice (and therefore the above union)   is typically infinite. For  instance, for $n=3$ 
   we start with  a finite set of lines in $\RR P^2$,  form all their intersection points, then draw new lines through these
   points and so on.

       \vskip .2cm
 
   (c)  Let $V_\RR = \RR^n$ with coordinates $x_1,\cdots, x_n$. Take $\Hc$ to be the arrangement
    of the following hyperplanes:
    \[
    \{x_i=0\}, \,\, i=1,\cdots, n,\quad \{x_i=x_{i+1}\}, \,\, i=1,\cdots, n-1. 
    \]
    There are ${n+1}\choose 2$ one-dimensional flats of $\Hc$, they have the form
    \[
    L_{[i,j]} \,=\,\bigl\{ x \,\bigl| \, x_i=x_{i+1}=\cdots = x_j;  \,\,\, x_k=0, k\notin [i,j]\bigr\},
    \quad 1\leq i\leq j\leq n. 
    \]
   On the other hand, consider $\RR^{n+1}$ with coordinates $y_0, \cdots, y_n$ and 
   let  $W_\RR = \RR^{n+1}/\RR\cdot (1,\cdots, 1)$.  Thus   $W_\RR = \hen^*$ is the space of
   weights for the Lie algebra $\sen\len_{n+1}(\RR)$. 
   We have an isomorphism $V_\RR\to W_\RR^*$ which takes the $i$th basis vector $e_i\in V_\RR$,
   $i=1,\cdots, n$, to the functional $y\mapsto y_{i-1}-y_i$ (simple co-root). 
   This isomorphism takes $L_{[i,j]}$ to the co-root hyperplane $\{y_{i-1}=y_j\}$.
   Therefore the dual arrangement $\Hc^\vee$ is the co-root arrangement in $\hen^*$. 
   
   \vskip .2cm
   
   Next,   flats of $\Hc^\vee$ are in bijection with equivalence relations
   $R$ on the set $\{0,1,\cdots, n\}$. The flat corresponding  to $R$ has the form
   \[
   M_R \,=\, \bigl\{ y\,\bigl| \, y_i=y_j \text{ whenever }i\equiv_R j.\bigr\}. 
   \]
 It is one-dimensional if and only if $R$ has only $2$ equivalence classes,
   both non-empty.
   Thus there are $2^{n-1}-1$ one-dimensional flats of $\Hc^\vee$ and so the double dual
   arrangement $\Hc^{\vee\vee}$ consists of $2^{n-1}-1$ hyperplanes and is much bigger than
   $\Hc$. 
   
   \vskip .2cm

   \end{exas}
    
    \begin{prop}\label{prop:FS-smooth}
    (a) If $\Gc\in D^b(V_\RR, \Sc_\RR)$, then  $\FS(\Gc) \in D^b(V^*_\RR, \Sc^\vee_\RR)$. 
    
    \vskip .2cm
    
    (b) If $\Fc\in \Perv(V_\CC, \Sc_\CC)$, then $\FS(\Fc)[-n]\in \Perv(V^*_\CC, \Sc^\vee_\CC)$. 
    \end{prop}
    
    \noindent{\sl Proof:}   As in the proof   of Proposition \ref{prop:spec-smooth}, 
the real and complex case are completely parallel, so we treat  the real case,
dropping the subscript $\RR$. The microsupport of $\Gc$ is contained in the union
of the $T^*_L V = L\times L^\perp$ over $L\in\Fl(\Hc)$. Now, the effect of $\FS$
on microsupports is  via the identification (``Legendre transform'')
\[
T^*V \,=\, V\times V^* \lra V^*\times V \,=\, T^* V^*.
\]
This identification takes $T^*_LV$ to  $T^*_{L^\perp} V^*$. This means that $\FS(\Gc)$ is smooth with respect
to the stratification $\Sc^*$ formed by the generic parts

 \[
L^\perp_\circ \,\, = \,\, L^\perp\,\,  \setminus\,\,  \bigcup_{L_1^\perp \not\supset L^\perp} L_1^\perp, 
\quad L\in\Fl (\Hc). 
\]
Now, $\Sc^\vee$  refines $\Sc^*$, so $\FS(\Gc)$ is smooth with respect to
$\Sc^\vee$. \qed

    \paragraph{C. Big and small dual cones.}
   Let $A^\vee\in \Sc^\vee_\RR$ be a  face. Its {\em big dual cone} is defined as  
    \be\label{eq:dual-cone}
     U(A^\vee) \,\,=\,\, \bigl\{ x\in V \,\bigl| f(x) \geq 0, \,\,\forall f\in A^\vee \bigr\} \,\,\subset \,\, V_\RR. 
     \ee
    It is a closed polyhedral cone in $V$ with nonempty interior, the union of the closures
     of (in general, several)  chambers  of $\Hc$. 
     
\begin{figure}[H]
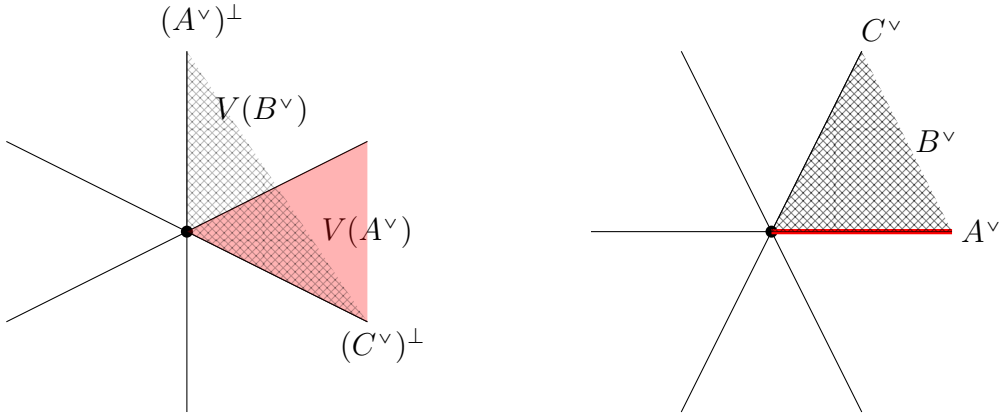

   \centering

  \btp[scale=.4]
  
  \node (0) at (0,0){};
  \fill[] (0) circle (0.2); 
  \draw (0,6) -- (0,-6); 
  \draw (6,3) -- (-6,-3);
  \draw (6,-3) -- (-6,3); 
  
  \node at (.5, 7) {$(A^\vee)^\perp$};
  \node at (6.5, -3.8){$(C^\vee)^\perp$}; 
\node at (6,0){$V(A^\vee)$};
\draw[fill=red, line width=0, opacity=0.3] (6,3) -- (0,0) -- (6,-3); 
\draw[fill=green, pattern=crosshatch, line width=0, opacity=0.5] (0,6) -- (0,0) -- (6,-3); 
\node at (2.5,4){$V(B^\vee)$};

  \etp
  \quad\quad\quad\quad\quad
   \btp[scale=.4]
  
  \node (0) at (0,0){};
  \fill[] (0) circle (0.2); 
  \draw (0,0) -- (-6,0); 
  \draw (3,6) -- (-3,-6);
  \draw (3,-6) -- (-3,6); 
  
  \node at (7,0) {$A^\vee$}; 
  \node at (3.7,6.7){$C^\vee$}; 
    \node at (5.5,3) {$B^\vee$}; 
    \draw[color=red, line width=2] (0,0) -- (6,0); 
    \draw[fill=green, pattern=crosshatch,  line width=0, opacity=0.7] (6,0) -- (0,0) -- (3,6);

  \etp
  
  \caption{Small dual cones.}\label{fig:small}
  \end{figure}

     The {\em small dual cone} of $A^\vee$ is defined as
     \be\label{eq:small-cone}
     V(A^\vee) \,\,=\,\, \bigcap_{B^\vee \geq A^\vee \text{ chamber}} U(B^\vee). 
     \ee
      It is a strictly convex  (not containing $\RR$-linear subspaces) closed polyhedral cone in $V_\RR$. 
      Note that $U(A^\vee)=V(A^\vee)$ if $A^\vee$ is a chamber but $U(A^\vee)$  can be strictly larger than
     $V(A^\vee)$ in general. For example, if $A^\vee$ is a half-line ($1$-dimensional face) of $\Hc^\vee$,
     then $U(A^\vee)$ is a closed half-space in $V_\RR$, while $V_A$ is strictly convex, cf. Fig.
     \ref{fig:small}.

  The next statement is clear from the definitions.
  
  \begin{prop}\label{prop:UV-monot}
   If $A^\vee_1 \leq A^\vee_2$, then $U(A^\vee_1) \supset U(A^\vee_2)$
   and  $V(A^\vee_1) \subset V(A^\vee_2)$. \qed
  \end{prop}
  
  \begin{prop}\label{prop:V(A)-pos}
   Let $f\in A^\vee$ be arbitrary. Then: 
   
   \vskip .2cm
   
   (a)  $U(A^\vee)$ is the union of all faces $B$ of $\Hc$ such that $f|_B\geq 0$ (non-strict inequality). 

   \vskip .2cm
          
  (b)    $V(A^\vee)$ is the union of $0$
     and all the faces $B$ of $\Hc$ such that $f|_B >0$ (strict inequality everywhere).

        \end{prop}
     
     \noindent{\sl Proof:} 
      (a)  Since $A^\vee$ is a face of $\Hc^\vee$, for each $f\in A^\vee$ the pattern of signs
       (positive, negative or zero) of $f$ on faces of $\Hc$ is the same. So the requirement that $f|_B\geq 0$
       for each $f\in A^\vee$ (appearing in the definition of $U(A^\vee)$) is equivalent 
       to the requirement that $f|_B\geq 0$ for any particular choice of $f\in A^\vee$ (appearing
       in the statement of the proposition). 
            
       \vskip .2cm
     (b) Let $V'$ be the union of the faces in question. If 
    $B\neq 0$ is a face of $\Hc$
     such that $B\subset V'$, i.e., that $f|_B>0$, 
     then $g|_B>0$ for any  $g\in A^\vee$,  by definition of the dual arrangement.
     This means that for any $B^\vee \geq A^\vee$ and any $g\in B^\vee$ sufficiently
     close to $A^\vee$ we still have $g|_B>0$. This, further implies (again, by the definition
     of the dual arrangement) that  
     for any $B^\vee \geq A^\vee$ and any $g\in B^\vee$ whatsoever 
     we still have $g|_B>0$. This means that  $B\subset V(B^\vee)$ for any $B^\vee \geq A^\vee$,
  in other words, that $B\subset V(A^\vee)$. We proved that $V'\subset V(A^\vee)$. 
  
  Conversely, suppose $B\subset V(A^\vee)$. For any chamber $B^\vee\geq A^\vee$
  and any $g\in B^\vee$ the restriction $g|_B$ cannot vanish, since that would mean
  that $g$ is not inside a chamber of a dual arrangement. Therefore $g|_B>0$ everywhere.
  Now, if $f\in A^\vee$ and $A^\vee$ is not a chamber, then looking at $g$ varying
  in a small transverse ball to $A^\vee$ near $f$ in $V^*_\RR$, we see that all such $g|_B$
  must be positive and therefore $f|_B$ must be positive. In other words,
  we proved that $V(A^\vee)\subset V'$. 
  
  \qed
  
       \begin{cor}\label{prop:UandV}
    We have
      \[
     U(A^\vee) \,\,=\,\,\bigcup_{B^\vee \geq A^\vee} U(B^\vee)
     \,\,=\,\, \bigcup_{B^\vee \geq A^\vee} V(B^\vee). 
     \]
         \qed
     \end{cor}
     
  We now analyze the nature of the covering of $U(A^\vee)$ by the $U(B^\vee)$, $B^\vee\geq A^\vee$.  
  All $B^\vee\geq A^\vee$ are in bijection with faces of the {\em quotient arrangement}
  $\Hc^\vee/A^\vee$ in the quotient space $V^*_\RR/\Lin_\RR(A^\vee)$, cf. \cite{KS} \S 2B. 
  We denote by $B^\vee/A^\vee$ the face of $\Hc^\vee/A^\vee$ corresponding to $B^\vee\geq A^\vee$. 
  
  \begin{prop}\label{prop:covering-U}
  Let  $A^\vee\in \Sc^\vee_\RR$ and $B\in \Sc_\RR$. Then:
  
  \vskip .2cm
  
  (a) There is a closed convex polyhedral cone $K(A^\vee, B)\subset V^*_\RR/\Lin_\RR(A^\vee)$,
  a union of faces of $\Hc^\vee/A^\vee$,  which has the following property:
  \[  
\text{  For $B^\vee\geq A^\vee$ we have
   $B\subset U(B^\vee)$ if and only $B^\vee/A^\vee\subset K(A^\vee, B)$.}
   \]
   
   \vskip .2cm
   
   (b) The cone $K(A^\vee, B)$ coincides with the whole $V^*_\RR/\Lin_\RR(A^\vee)$
   if and only if $B\subset V(A^\vee)$. 
  
  \end{prop} 
  
  \noindent {\sl Proof:} (a) Let $U(B)\subset V^*_\RR$ be the dual cone to $B$, i.e., the
  set of $f\in V^*_\RR$ such that $f|_B\geq 0$. It is a convex, closed polyhedral cone in
  $V^*_\RR$ which is a union of faces of $\Hc^\vee$. 
   In fact, the condition  $B^\vee \subset U(B)$
 is equivalent to  $B\subset U(B^\vee)$, 
 both meaning that $(b^\vee, b)\geq 0$ for each
  $b^\vee\in B^\vee$ and $b\in B$.  
  
  \vskip .1cm
  
  Let also $(V^*_\RR)^{\geq A^\vee}\subset V^*_\RR$ be the union of all faces $B^\vee$
  of $\Hc^\vee$ such that $B^\vee\geq A^\vee$. It is a convex, open polyhedral cone in $V^*_\RR$.
  The intersection $U(B)\cap (V^*_\RR)^{\geq A^\vee}$ is then a convex polyhedral
  cone which is closed in $(V^*_\RR)^{\geq A^\vee}$. Since this cone is a union of faces
  $B^\vee\geq A^\vee$, it projects to a convex closed polyhedral cone in
  $V^*_\RR/\Lin_\RR(A^\vee)$ which we denote $K(A^\vee, B)$. By construction,
  $K(A^\vee, B)$ satisfies the required property.
  
  \vskip .2cm
  
  (b) This is a reformulation of the formula \eqref{eq:small-cone} defining $V(A^\vee)$. \qed
  
  \vskip .2cm
  
  Note an appealing numerical corollary of Proposition \ref{prop:covering-U}. For any subset $Z\subset V$ we
  denote by $\1_Z: V\to\RR$ its characteristic function, equal to $1$ on $Z$ and to $0$ elsewhere. 
  
  \begin{cor}\label{cor:alternating} (inclusion - exclusion formulas)
  We have the identities
  \[
  \1_{U(A^\vee)} \,\,=\,\,\sum_{B^\vee\supset A^\vee} (-1)^{\dim(B^\vee)-\dim(A^\vee)} \1_{V(A^\vee)}, 
  \leqno (a)
  \]
  \[
  \1_{V(A^\vee)} \,\,=\,\,\sum_{B^\vee\supset A^\vee} (-1)^{\dim(B^\vee)-\dim(A^\vee)} \1_{U(A^\vee)}. 
  \leqno (b)
 \]
  \end{cor}
  Identities of this general nature (representing the characteristic function of a convex polytope as
  an alternating sum of characteristic functions of simplices or cones) are familiar in the theory of
  convex polytopes \cite{V}  \cite{FL} and the theory of automorphic forms, see, e.g.,  \cite{Ar}, \S 11. 
  
  We note the similarity of the identities (a) and (b) with Proposition \ref{prop:stalks and hyperstalks}(b) and 
  Corollary \ref{cor:stalks-and-hyper} relating the usual stalks and hyperbolic stalks of a perverse sheaf.
  In fact, we will use a ``categorified'' version of these identities to relate the usual and hyperbolic stalks
  of the Fourier-Sato transform.

  \vskip .2cm
  
  \noindent {\sl Proof of Corollary \ref{cor:alternating} :} (b) Write the RHS of the proposed identity as
  $\sum_B c_B \1_B$ with $B$ running over faces of $\Hc$.  Part (a) of Proposition   \ref{prop:covering-U}
  implies that 
  \[
  c_B \,  = \sum_{B^\vee \geq A^\vee\atop B^\vee/A^\vee \subset K(A^\vee, B)} (-1)^{\dim(B^\vee) -\dim(A^\vee)}\,\,=\,\,
  (-1)^{n-\dim(A^\vee)} \chi \bigl(H^\bullet_c(K(A^\vee, B), \k\bigr)
  \] 
  is the signed (calculated from the top)
   Euler characteristic of the cohomology with compact support of  the cone $K(A^\vee, B)$.
   This signed Euler characteristic is equal to $0$ unless $K(A^\vee, B)$ is the entire vector space, in which case it is
 $1$. By Part (b) of Proposition     \ref{prop:covering-U}  this happens precisely when
  $B\subset V(A^\vee)$,  so  the identity is proved. 
  
  \vskip .2cm
  
  (a) is a formal consequence of (b) in virtue of the identity
   \[
  \sum_{B^\vee \geq A^\vee} (-1)^{\dim(B^\vee)-\dim(A^\vee)}\,\,=\,\, 1
  \]
  (the Euler characteristic of the link of $A^\vee$). \qed

     \paragraph{D. The real result.}
 
    \begin{thm}\label{thm:FS-real}
     Let $\Gc\in D^b(V_\RR, \Sc_\RR)$ be represented by  a complex $(\Gc_A, \gamma_{AB})$
     of representations of $(\Sc_\RR, \leq)$.

     \vskip .2cm
     
     (a)  The stalk $\FS(\Gc)_{A^\vee}$ of  $\FS(\Gc)$ at a face
      $A^\vee\in\Sc^\vee_\RR$ is identified with the complex
     \[
     \Uc_{A^\vee} \,\, :=  \,\, \Tot\,\biggl\{ \bigoplus_{\dim(B)=0, \atop B\subset U(A^\vee)} \Gc_B \otimes
     \OR(B) \buildrel \gamma\otimes\eps \over\lra
     \bigoplus_{\dim(B)=1, \atop B\subset U(A^\vee)} \Gc_B \otimes
     \OR(B) \buildrel \gamma\otimes\eps \over\lra
\cdots\biggr\}. 
     \]
     
      \vskip .2cm
     
     (b) Let $A^\vee_1\leq A^\vee_2$ be two faces of $\Hc^\vee$.
     Then the inclusion 
      $U(A^\vee_1)\supset
     U(A^\vee_2)$ (Proposition \ref{prop:UV-monot}) exhibits $\Uc(A^\vee_2)$
     as a quotient complex of $\Uc(A^\vee_1)$, and the 
     the generalization map $\gamma_{A^\vee_1,  A^\vee_2}: \FS(\Gc)_{A^\vee_1}
     \to \FS(\Gc)_{A^\vee_2}$ 
     of $\FS(\Gc)$ 
     is identified with the quotient map $\Uc_{A_1^\vee}\to \Uc_{A_2^\vee}$.      
        \end{thm}
     
     \noindent {\sl Proof:} (a) 
     Let $f\in  A^\vee$.  By Proposition \ref{prop:FS-stalk-K} we have
 
      \[
    \FS(\Gc)_{A^\vee}\,\simeq \,   \FS(\Gc)_f  \,\simeq \, R\Gamma_c(P_f,\Gc). 
      \] 
      We now use the resolution of $\Gc$ given by Proposition \ref{prop:LA-res}(ii). The $p$th term
      of this resolution is the direct sum of $j_{B!}\,  \ul{\Gc_B}_B[p]$ where $B$ runs over 
      $p$-dimensional faces of $\Hc$. 
      
      \begin{lem}
      Let $B$ be a face of $\Hc$ and $E$ be any $\k$-vector space. We have
      natural quasi-isomorphisms
      \[
   R\Gamma_c(P_f, j_{B!} \ul E_B)[\dim B]  \simeq    \begin{cases}
    E\otimes \orr(B), & \text{if } B\subset P_f, 
    \\
    0, & \text{if } B\not\subset P_f. 
      \end{cases}
      \]
      \end{lem}     
      
      \noindent {\sl Proof of the lemma:} The case $B\subset P_f$ follows from the  canonical identification
       $R\Gamma_c(B, \k) \simeq \orr(B)
      [-\dim(B)]$
      (compactly supported cohomology of a cell with constant coefficients). 
      Suppose $B\not\subset P_f$. If $B$ does not meet $P_f$ at all, then the statement is obvious.
      If $B$ does meet $P_f$, then the intersection $B\cap P_f$ is homeomorphic to a closed half-space
      in a Euclidean space, i.e., to a Cartesian product of several open intervals $(0,1)$ and one
      half-open interval $[0,1)$. So our statement follows from the fact that
      $H^\bullet_c( [0,1), \k) = 0$. \qed

  \vskip .2cm
  
  Applying this lemma to the resolution of  $\Gc$ given by Proposition \ref{prop:LA-res}(ii), we obtain
  a complex representing $R\Gamma_c(P_f, \Gc)$ whose $p$th term is the sum of the $\Gc_B\otimes \orr_B$ 
  for $B$  running over
  $p$-dimensional faces $B\subset P_f$ and the differential is formed by the maps $\gamma\otimes\eps$. 
  By Proposition \ref{prop:covering-U} (a), the condition $B\subset P_f$
  is equivalent to $B\subset U(A^\vee)$.  This proves part (a) of Theorem \ref {thm:FS-real}. 
  
  \vskip .2cm
  
  We now prove part (b). Let $f_1\in A_1^\vee$ and $f_2\in A_2^\vee$ be a small deformation of $f_1$. 
  As in the proof of (a), we can write our generalization map as
  \[
  \gamma_{A^\vee_1, A^\vee_2}: R\Gamma_c(P_{f_1}, \Gc) \lra  R\Gamma_c(P_{f_2}, \Gc). 
  \]
  As before,  consider first the case $\Gc = j_{B!} \ul{ E}_B$ for some face $B$
  and some $\k$-vector space $E$. In this case we find that $\gamma_{A^\vee_1, A^\vee_2}$
  is equal to the identity map, 
  if $B$ is contained in $P_{f_2}$ (and therefore in $P_{f_1}$),  and it is equal to $0$ otherwise
  (since the target is the zero vector space). 
  That is, claim (b) obviously holds in this case. The case of general $\Gc$ is now  obtained
  from this by considering the resolution  of $\Gc$ given by Proposition \ref{prop:LA-res}(ii). 
  Theorem \ref {thm:FS-real} is proved.

     \paragraph{E. The complex result.} 
     
 Let $\Fc\in\Perv(V_\CC, \Sc_\CC)$ correspond to a hyperbolic sheaf 
     $\Qc=(E_A, \gamma_{AB}, \delta_{BA})$. By Proposition \ref {prop:FS-smooth},
      $\FS(\Fc)[-n]$ lies in $\Perv(V_\CC^*, \Sc^\vee_\CC)$ and so is described by
      a  hyperbolic sheaf which we denote
      $\Qc^\vee = (E^\vee_{A^\vee}, \gamma_{A^\vee , \Ap^\vee}, \delta_{\Ap^\vee,  A^\vee})$. 
      Here $A^\vee \leq \Ap^\vee$ are faces of the arrangement $\Hc^\vee$. 
      
      It turns out that the hyperbolic stalks $E^\vee_{A^\vee}$ are governed by the small dual cones
      $V(A^\vee)$.

     \begin{thm}\label{thm:hyp-stalk-FS}
(a) The space $E^\vee_{A^\vee}$ is quasi-isomorphic to the complex
\[
\Vc_{A^\vee} \,= \,\,\biggl\{\bigoplus_{\dim(B)=0, \atop B\subset V (A^\vee)} E_B \otimes
     \OR_{V/B} \buildrel \gamma\otimes\eps \over\lra
     \bigoplus_{\dim(B)=1, \atop B\subset V(A^\vee)} E_B \otimes
     \OR_{V/B} \buildrel \gamma\otimes\eps \over\lra
\cdots\biggr\}. 
\]
In other words,  $\Vc_{A^\vee}$ is exact everywhere except the leftmost term, where the
cohomology (kernel) is identified with $E^\vee_{A^\vee}$.

\vskip .2cm

(a') The space $E^\vee_{A^\vee}$ is also quasi-isomorphic to the complex
\[
\Vc_{A^\vee}^\dagger  \,= \,\,\biggl\{\cdots  \buildrel \delta\otimes\eps \over\lra
\bigoplus_{\dim(B)=1, \atop B\subset V(A^\vee)} E_B \otimes
     \OR_{V/B} 
      \buildrel \delta\otimes\eps \over\lra
\bigoplus_{\dim(B)=0, \atop B\subset V (A^\vee)} E_B \otimes
     \OR_{V/B} \biggr\}. 
\]
 In other words,  $\Vc_{A^\vee}^\dagger$
  is exact everywhere except the rightmost term, where the
cohomology (cokernel) is identified with $E^\vee_{A^\vee}$.

\vskip .2cm

(b) Let $A^\vee_1\leq  A^\vee_2$ be two faces of $\Hc^\vee$. Then the embedding
$V(A^\vee_1) \subset V(A^\vee_2)$ realizes $\Vc_{A^\vee_1}$ as a quotient complex of
$\Vc_{A^\vee_2}$, and the map $\delta_{A^\vee_2, A^\vee_1}: E^\vee_{A_2^\vee}
\to E^\vee_{A^\vee_1}$ is identified with the quotient map $\Vc_{A_2^\vee}\to
\Vc_{A_1^\vee}$. 

\vskip .2cm

(b') In the situation of (b), the embedding
$V(A^\vee_1) \subset V(A_2^\vee)$ realizes $\Vc^\dagger_{A^\vee_1}$ as a subcomplex of
$\Vc_{A_2^\vee}^\dagger$, and the map $\gamma_{A^\vee_1, A_2^\vee}$
is identified with the embedding $\Vc^\dagger_{A^\vee_1} \to \Vc^\dagger_{A_2^\vee}$. 
 
   \end{thm}
   
   \begin{rems}
  (a) Note that for $A^\vee=0$, the cone $V(A^\vee)$ is equal to $\{0\}$, therefore
  $E^\vee_0$ is identified with $E_0$.  
  
  \vskip .2cm

  (b) Let $A^\vee\neq 0$. Then, by Proposition \ref{prop:V(A)-pos}(b) one can re-write the complex 
  $\Vc_{A^\vee}$ as
    \[
  E_0\otimes\orr_V \buildrel \gamma\otimes\eps\over\lra \bigoplus_{\dim(B)=1\atop f|_B>0} 
  E_B\otimes \orr_{V/B}  \buildrel \gamma\otimes\eps\over\lra  
  \bigoplus_{\dim(B)=2\atop f|_B>0} 
  E_B\otimes \orr_{V/B}  \buildrel \gamma\otimes\eps\over\lra \cdots, 
  \]
  where $f\in A^\vee$ is an arbitrary element. Similarly for $\Vc^\dagger_{A^\vee}$. 
   \end{rems}
   
   The proof  of Theorem \ref{thm:hyp-stalk-FS} is based on the following preliminary result which shows that  the big dual cones
   $U(A^\vee)$ govern the ordinary stalks, not hyperbolic stalks of $\FS(\Fc)$, 
   
   \begin{prop}\label{prop:stalks-FS-comp}
   (a) If $A^\vee\in\Sc_\RR$ is any face, then the ordinary stalk $\FS(\Fc)_{A^\vee}$ is quasi-isomorphic to the complex
   \[
   \FS(\Fc)_{A^\vee} \,\,\simeq\,\, \biggl\{ \bigoplus_{\dim(B)=0, \atop B\subset U(A^\vee)} E_B \otimes
     \OR_{V/B} \buildrel \gamma\otimes\eps \over\lra
     \bigoplus_{\dim(B)=1, \atop B\subset U(A^\vee)} E_B \otimes
     \OR_{V/B} \buildrel \gamma\otimes\eps \over\lra
\cdots
\biggr\}.
 \]
   
 \vskip .2cm
 
 (b) The generalization maps for the $\FS(\Fc)_{A^\vee}$ are induced by
 the projections of the complexes in (a), similarly to Theorem \ref{thm:FS-real}(b). 
   \end{prop}
   
   \noindent{\sl Proof of Proposition \ref {prop:stalks-FS-comp}:} 
    Our statement   will follow from Theorem \ref{thm:FS-real}, if we establish the following.
   
   \begin{prop}\label{prop:FS-i!}
   For  any $\Fc\in\Perv(V_\CC, \Sc_\CC)$ we have an identification
   \[
   \FS(i_\RR^!\Fc) \,\,\simeq \,\, i_{\RR}^*\,  \FS (\Fc),
   \]
   where $i_\RR$ on the right means the embedding $V_\RR^* \to V_\CC^*$. 
   \end{prop}
     
    \noindent {\sl Proof of Proposition \ref{prop:FS-i!}:}  We first recall the behavior of the Fourier-Sato
     transform
    with respect to an arbitrary $\RR$-linear map $\phi: W_1\to W_2$  of $\RR$-vector spaces.
    Denoting $^t\phi: W_2^*\to W_1^*$ the transposed map, we have, for any conic complex $\Gc$ on $W_2$:
    \[
    \FS(\phi^! \Gc) \,\,\simeq \,\, R(^t\phi)_* \, \FS(\Gc)
    \]
    see \cite{KaSha} Prop. 3.7.14. 
    
    We specialize this to $\phi=i_\RR: V_\RR\to V_\CC$ and $\Gc=\Fc$. 
      In this case 
      \[
      ^t \phi \,\,= \,\, \Re: V^*_\CC\lra V^*_\RR
      \]
      is the real part map. So after replacing $V^*$ by $V$ and $\FS(\Fc)$ by $\Fc$,
       Proposition \ref{prop:FS-i!} reduces to the following.
      
      \begin{lem}\label{prop:i!-RRe}
      For  any $\Fc\in D^b(V_\CC, \Sc_\CC)$ we have an identification 
     $ i_\RR^*\Fc \,\simeq \, R \, \Re_*(\Fc) $,  
      where $\Re: V_\CC\to V_\RR$ is the real part map for $V$. 
      \end{lem}
      
      \noindent {\sl Proof of Lemma \ref{prop:i!-RRe}:} 
    We consider $\Re: V_\CC\to V_\RR$
      as a real vector bundle over $V_\RR$. The complex $\Fc$, being constructible with
      respect to the complexification of a real hyperplane arrangement,   is conic with respect to this vector bundle structure.
      Therefore the stalk at $x\in V_\RR$ of $i_\RR^*\Fc$ which is $R\Gamma(U, \Fc)$
      for a small open $U\subset V_\CC$ containing $x$, is equal to $R\Gamma(\Re^{-1}(U\cap V_\RR), \Fc)$
      which is the stalk of $R\, \Re_*(\Fc)$ at $x$. \qed 
      
      This finishes the proof of Propositions \ref{prop:FS-i!} and \ref{prop:stalks-FS-comp}. 
      
      \vskip .2cm
      
      \noindent{\sl Proof of Theorem \ref{thm:hyp-stalk-FS}:} We prove (a') and (b'). Parts (a)
      and (b) follow by Verdier duality. 
      
      We denote $\Kc=\FS(\Fc)$, and let  $\Lc = \Kc^* = \FS(\Fc^*)$ be the Verdier dual
      perverse sheaf. By definition, $E^\vee_{A^\vee}$ is the stalk at $A^\vee$ of
      \[
      i_\RR^!\Kc\,\simeq \, (i_\RR^*\Lc)^*\, \simeq \, \bigl(i_\RR^* \FS(\Fc^*)\bigr)^*. 
      \]
      First, we recall that
       $\Fc^*$ is represented by the hyperbolic sheaf $(E_A^*, \delta^*_{BA}, \gamma_{AB}^*)$. 
      Applying Proposition \ref{prop:stalks-FS-comp} to $\Fc^*$ we write the 
     stalk of $i_\RR^*\Lc$ at $A^\vee$ as
     \be\label{eq:complex-L} 
     \Lc_{A^\vee} \,\,\simeq \,\,\biggl\{ \bigoplus_{\dim(B)=0\atop B\subset U(A^\vee)} 
     E_B^* \otimes \OR_{V/B} \buildrel\delta^*\otimes\eps \over \lra 
     \bigoplus_{\dim(B)=1\atop B\subset U(A^\vee)} 
     E_B^* \otimes \OR_{V/B}  \buildrel\delta^*\otimes\eps \over \lra \cdots
     \biggr\}. 
     \ee
     Further,
      for $A^\vee_1\leq A^\vee_2$ we have $U(A^\vee_1)\supset U(A^\vee_2)$ and 
       Proposition \ref{prop:stalks-FS-comp} implies that the 
     generalization map $\digamma_{A^\vee_1, A^\vee_2} : \Lc_{A^\vee_1}\to\Lc_{A^\vee_2}$
     is given by the projection of the corresponding complexes in \eqref{eq:complex-L}.
     
     We now recall the following general procedure on finding the
     stalks and generalization maps of the Verdier dual complex. See, e.g., \cite{KS} Prop. 1.11.
     We formulate it here for complexes on $V^*_\RR$  constructible with respect to $\Sc^\vee_\RR$. 
     
     \begin{lem}\label{lem:verdier}
     Let $\Mc\in D^b(V_\RR^*, \Sc^\vee_\RR)$ correspond to a complex
     $(\Mc_{A^\vee}, \digamma_{A_1^\vee, A^\vee_2})$ of representations of
     $\Sc^\vee_\RR$. Then:
     
     \vskip .2cm
     
     (a) The stalk of $\Mc^*$ at $A^\vee$ is identified with the complex
     \[
     \Dc_{A^\vee} \,\,=\,\, \Tot\biggl\{ \cdots \buildrel \digamma^*\otimes\eps^* \over\lra
      \bigoplus_{C^\vee >_1 A^\vee} (\Mc_{C^\vee})^* \otimes\OR_{C^\vee} 
      \buildrel \digamma^*\otimes\eps^* \over\lra
       (\Mc_{A^\vee})^* \otimes\OR_{A^\vee}\biggr\},
       \]
       with the horizontal grading associating to the summand $  (\Mc_{A^\vee})^* \otimes\OR_{A^\vee}$
       degree $-\dim(A^\vee)$.  The  horizontal differential $\digamma^*\otimes\eps^*$
      has, as the matrix element corresponding  to $C_2^\vee >_1 C_1^\vee \geq A^\vee$, 
       the tensor product of the dual maps to $\digamma_{C^\vee_1, C^\vee_2}$
      and to  $\eps_{C^\vee_1, C^\vee_2}$.
       
       \vskip .2cm
       
       (b) For two faces $A^\vee_1\leq A^\vee_2$ the generalization map 
       $(\Mc^*)_{A^\vee_1} \to (\Mc^*)_{A^\vee_2}$ of $\Mc^*$, is identified with the projection
       of the complexes $\Dc_{A^\vee_1} \to\Dc_{A^\vee_2}$. \qed
     \end{lem}

      \vskip.2cm
      Applying part (a) of  the lemma to $\Mc=\Lc$ and
      substituting, instead of each $\Lc_{C^\vee}$, its expansion \eqref{eq:complex-L},
      we identify (quasi-isomorphically) $E^\vee_{A^\vee}$   with
       the total complex of the following double complex. We denote this total complex
       $\Ec_{A^\vee}$.
     \[
      \xymatrix{
      \cdots \ar[r]^{\hskip -3.5cm \emb\otimes\eps}& \bigoplus_{C^\vee >_1 A^\vee} \bigoplus_{\dim(B)=0\atop B\subset U(C^\vee)} 
      E_B \otimes\OR_{V/B} \otimes\OR_{C^\vee} \ar[r]^{\hskip .5cm \emb\otimes\eps}& 
     \bigoplus_{\dim(B)=0\atop B\subset U(C^\vee)} 
      E_B \otimes\OR_{V/B} \otimes\OR_{A^\vee} 
      \\
        \cdots \ar[r]^{\hskip -3.5cm \emb\otimes\eps}
        & \bigoplus_{C^\vee >_1 A^\vee} \bigoplus_{\dim(B)=1\atop B\subset U(C^\vee)} 
      E_B \otimes\OR_{V/B} \otimes\OR_{C^\vee}
         \ar[u]_{\delta\otimes\eps\otimes\Id}
       \ar[r]^{\hskip .5cm \emb\otimes\eps}& 
     \bigoplus_{\dim(B)=1\atop B\subset U(C^\vee)} 
      E_B \otimes\OR_{V/B} \otimes\OR_{A^\vee} 
      \ar[u]_{\delta\otimes\eps\otimes\Id}
      \\
        \cdots \ar[r]^{\hskip -3.5cm \emb\otimes\eps}
        & \bigoplus_{C^\vee >_1 A^\vee} \bigoplus_{\dim(B)=2\atop B\subset U(C^\vee)} 
      E_B \otimes\OR_{V/B} \otimes\OR_{C^\vee}
         \ar[u]_{\delta\otimes\eps\otimes\Id}
       \ar[r]^{\hskip .5cm \emb\otimes\eps}& 
     \bigoplus_{\dim(B)=2\atop B\subset U(C^\vee)} 
      E_B \otimes\OR_{V/B} \otimes\OR_{A^\vee} 
         \ar[u]_{\delta\otimes\eps\otimes\Id}
         \\
         & \vdots 
           \ar[u]_{\delta\otimes\eps\otimes\Id}
         & \vdots 
           \ar[u]_{\delta\otimes\eps\otimes\Id}
    }  
      \]
  Here the vertical differentials are dual to those in $\Lc_{C^\vee}$, i.e., given
  by the $\delta$ maps.    Matrix elements of the horizontal differential
   are dual to the $\digamma$ maps for
  $\Lc$, and those $\digamma$ maps are given by the projections. So each matrix element in
  question is in fact the product of an {\em embedding}  of $\delta$-complexes and the $\eps$
  map of orientation torsors.
  
  For two faces $A^\vee_1 \leq A^\vee_2$ the generalization map
  $\gamma_{A^\vee_1, A^\vee_2}:  E^\vee_{A^\vee_1}\to E^\vee_{A^\vee_2}$ is identified,
  by part (b) of Lemma \ref{lem:verdier},  with  the  projection
    $\Ec_{A^\vee_1}\to  \Ec_{A^\vee_2}$. 
  
  We now compare  $\Ec_{A^\vee}$ with the
  complex $\Vc^\dagger_{A^\vee}$ from the formulation of Theorem 
  \ref{thm:hyp-stalk-FS}(a').   Let $B$ be a face of $\Hc$.
  The summand corresponding to $B$ in  $\Vc^\dagger_{A^\vee}$,
   is  either $E_B\otimes \OR_{V/B}$ or $0$ depending on whether $B\subset V(A^\vee)$
   or not. On the other hand, $\Ec_{A^\vee}$ has many summands associated to $B$,
   they are labelled by $C^\vee> A^\vee$ such that $B\subset U(C^\vee)$.
   By Proposition \ref{prop:covering-U}, such $C^\vee$ are in bijection with faces of the
  closed  polyhedral cone $K(A^\vee, B)$. So in the double complex above
   the summand $E_B\otimes \OR_{V/B}$ is multiplied by a combinatorial complex
   which is easily found to calculate the cohomology with compact support 
   $H^\bullet_c(K( A^\vee, B), \k)$. If $B\not\subset V(A^\vee)$, then,
   by the same Proposition  \ref{prop:covering-U},
    $K(A^\vee, B)$ is a proper closed cone with nonempty interior in
   $V^*_\RR/\Lin_\RR(A^\vee)$ and so its cohomology with compact support
   vanishes entirely. If $B\subset V(A^\vee)$, then
   $K(A^\vee, B) =  V^*_\RR/\Lin_\RR(A^\vee)$ so it has the top cohomology with
  compact support identified with $\OR_{V/A^\vee}$, so  the part of  $\Ec_{A^\vee}$
  corresponding to $B$ is quasi-isomorphic to $E_B\otimes \OR_{V/B}$. Moreover,
  we see that these quasi-isomorphisms combine into a quasi-isomorphism between
 $\Ec_{A^\vee}$  and $\Vc^\dagger_{A^\vee}$. This shows part (a') of Theorem 
  \ref{prop:stalks-FS-comp}. Part (b') follows by noticing that the projection
  $\Ec_{A^\vee_1}\to \Ec_{A^\vee_2}$ corresponds, under our quasi-isomorphism, to
  the embedding $\Vc^\dagger_{A^\vee_1}\to \Vc^\dagger_{A^\vee_2}$. 
  Theorem  \ref{prop:stalks-FS-comp} is proved.

   \vfill\eject
  
  \section{Applications to second microlocalization}

     \paragraph{ A.   Microlocalization. }  
     
       If $M\subset X$ is a $C^\infty$ submanifold of a $C^\infty$ manifold, as in \S \S \ref{sec:spec}A, then
      for any $\Gc\in D^b (X)$ the {\em microlocalization} of $\Gc$ along $M$ is defined as 
      \[
      \mu_M(\Gc) \,=\, \FS_M(\nu_M(\Gc)) \,\in\, D^b (T^*_MX), 
      \]
      see  \cite{KaSha} Ch. 4. 
     Here $\FS_M$ is the relative Fourier-Sato transform on the vector bundle $T_MX\to M$. 
     
       \vskip .2cm
     
     If $X=V$ is a real vector space with an arrangement $\Hc$, if  $\Gc\in D^b(V, \Sc_\RR)$
     and  $M$ is a vector subspace, then
     our descriptions of the Fourier-Sato transform and the specialization functors      
       can be combined to obtain a combinatorial description of  $\mu_M(\Gc)$. 
       We leave  this to the reader,  establishing instead some compatibility
      properties of various approaches to ``second microlocalization"
      of Kashiwara and Laurent, see \cite{laurent} and references therein.
For convenience we give a brief general introduction.

\paragraph{B. Iterated microlocalization.} 

\begin{lem}
Let $(W,\omega)$ be a symplectic $\RR$-vector space, and $L_1, L_2\subset W$ be Lagrangian
vector subspaces. Then the restriction of $\omega$ gives an identification
\[
\left( {L_1\over L_1\cap L_2}\right)^* \,\,\simeq \,\, \left( {L_2\over L_1\cap L_2}\right).
\]
\end{lem}

\noindent{\sl Proof:} Consider the restriction of $\omega$ to the subspace $L_1 + L_2$. Its kernel 
on this subspace is 
\[
(L_1+L_2)^\perp \,=\, L_1^\perp  \cap L_2^\perp \,=\, L_1\cap L_2. 
\]   
Therefore the restriction of $\omega$ makes
\[
{L_1+ L_2 \over L_1\cap L_2} \,\,=\,\, {L_1\over L_1\cap L_2} \,\,\oplus {L_2\over L_1\cap L_2}
\]    
into a symplectic vector space decomposed into the direct sum of two Lagrangian subspaces.
So these Lagrangan subspaces become dual to each other. \qed

\vskip .2cm

Let now $(S,\omega)$ be a $C^\infty$ symplectic manifold and $\Lambda_1, \Lambda_2\subset S$ be two
(smooth) Lagrangian submanifolds. We say that $\Lambda_1$ and $\Lambda_2$ 
{\em intersect cleanly} (in the symplectic sense),
 if,  locally near each  $x\in\Lambda_1\cap\Lambda_2$,    there is a
symplectomorphism of a neighborhood of $x$ in $S$  to a neighborhood of $0$ in
a symplectic vector space $W$, sending $\Lambda_i$ to linear Lagrangian subspaces $L_i$
as above. This implies that $\Lambda_1\cap\Lambda_2$ is smooth. 

\begin{cor}\label{cor:tt*}
If $\Lambda_1, \Lambda_2$ intersect cleanly, then the restriction of $\omega$ gives an identification
\[
T^*_{\Lambda_1\cap \Lambda_2} \Lambda_1 \,\,\simeq \,\, T_{\Lambda_1\cap \Lambda_2} \Lambda_2. \qed
\]
\end{cor}

Now let $X$ be a $C^\infty$ manifolds and $M,N\subset X$ be two smooth submanifolds. We assume that
they intersect cleanly in the sense that they can locally be brought by a diffeomorphism to two vector
subspaces in a vector space. Then $S=T^*X$ has two Lagrangian submanifolds
 $\Lambda_1= T_M^*X$, $\Lambda_2=T^*_NX$ which intersect cleanly in the symplectic sense. 
 Given a complex of sheaves $\Gc\in D^b(X)$, we have microlocalizations
 \[
 \mu_M(\Gc) \,\in \, D^b(\Lambda_1), \quad \mu_N(\Gc)\in D^b(\Lambda_2)
 \]
 and we can specialize and microlocalize further, getting two complexes of sheaves
 \[
 \mu_{\Lambda_1\cap \Lambda_2} \mu_M(\Gc) \,\,\in\,\, D^b(T^*_{\Lambda_1\cap \Lambda_2} \Lambda_1), 
 \quad 
 \nu_{\Lambda_1\cap \Lambda_2} \mu_N(\Gc) \,\,\in\,\, D^b(T_{\Lambda_1\cap \Lambda_2} \Lambda_2)
 \]
 on two spaces which are identified by Corollary \ref{cor:tt*}, so we can consider them
 as living on the same space. 
 One can then formulate
 
 \begin{secmic}
 Under which conditions on $M,N$ and $\Gc$  can we guarantee that
  \[
   \mu_{\Lambda_1\cap \Lambda_2}\, \mu_M(\Gc) \,\,\simeq \,\,  \nu_{\Lambda_1\cap \Lambda_2}\, \mu_N(\Gc) 
   \quad ?
  \]
 \end{secmic}

 \paragraph{C. Bi-microlocalization.} 

Let us restrict to the case $N\subset M$.  In this case we have

\begin{prop}
We have identifications
\[
T^*_{\Lambda_1\cap \Lambda_2} \Lambda_1 \,\,=  \,\, T_{\Lambda_1\cap \Lambda_2} \Lambda_2
\,\,\simeq \,\, T^*_NM \oplus (T^*_MX)|_N. 
\]
\end{prop}

\noindent{\sl Proof:} Obviously, $\Lambda_1\cap\Lambda_2$ projects, under $T^*X\to X$, to $N$.
Looking at the fibers of this projection, we find that $\Lambda_1\cap\Lambda_2 = (T^*_MX)|_N$. 
Looking at the Cartesian square
\[
\xymatrix{
\Lambda_1\cap\Lambda_2 \ar[r] 
\ar[d]_\rho
& T^*_MX=\Lambda_1
\ar[d]^\pi
\\
N \ar[r]& M
}
\]
with $\pi$ being a smooth fibration (projection of a vector bundle), we find that
\[
T^*_{\Lambda_1\cap\Lambda_2}\Lambda_1 \,\,\simeq \,\,\rho^* T^*_NM \,\,\simeq \,\,
T^*_NM \oplus (T^*_MX)|_N. \qed
\]

We already considered the situation of  a flag $N\subset M\subset X$  in discussing
bi-specialization $\nu_{NM}(\Gc)$ in
     \S \ref{sec:spec}E. 
Further, in this context Schapira and Takeuchi \cite{schapira-takeuchi} \cite{takeuchi}
     have defined the {\em bimicrolocalization} 
     \[
     \mu_{NM}(\Gc) \,=\, \FS_N (\nu_{NM}(\Gc)) \, \in \,  D^b(T^*_NM \oplus (T^*_MX)|_N). 
     \]
   Here $\FS_N$ is the relative Fourier-Sato transform on the
     vector bundle $T_NM \oplus (T_MX)|_N\to N$. So we have the following specialization-microlocalization diagram:
\be\label{eq:bimicro}
\xymatrix{
D^b(X) \ar[rr]^{\mu_M} 
\ar[d]_{\mu_N}
\ar[rrd]^{\mu_{NM}}
&& D^b(T^*_MX = \Lambda_1) 
\ar[d]^{\mu_{\Lambda_1\cap\Lambda_2}}
\\
D^b(T^*_NX = \Lambda_2) \ar[rr] _{\hskip -2.5cm \nu_{\Lambda_1\cap\Lambda_2}}
&& D^b\bigl( (T^*_NX\oplus (T^*_MX)|_N)) = T^*_{\Lambda_1\cap\Lambda_2}\Lambda_1 =
T_{\Lambda_1\cap\Lambda_2}\Lambda_2\bigr)
}
\ee
which gives  three possible ``second microlocalizations''. 

\paragraph{D. Comparisons  in the linear case. }

\begin{thm}\label{thm:micro-3}
Let $X=V$ be an $\RR$-vector space,
$N\subset M\subset V$ be  vector subspaces and   $\Hc$ an arrangement of hyperplanes
in $V$ with the corresponding face stratification $\Sc_\RR$. Then the diagram \eqref {eq:bimicro}
is canonically 2-commutative if we replace $D^b(V)$ with $D^b(V,\Sc_\RR)$. 
\end{thm} 

In the complex situation, when $V=V_\CC$ is a $\CC$-vector space, $N\subset\ M\subset V_\CC$ are
$\CC$-subspaces and $D^b(V,\Sc_\RR)$ is replaced by $\Perv(V_\CC, \Sc_\CC)$,
the commutativity of the outer square of  \eqref {eq:bimicro} was proved in \cite{FS}
using the $\Dc$-module techiques.

We will deduce Theorem \ref{thm:micro-3}
from the following result.

\begin{thm}[P. Schapira] \label{refo:micro-2}
Let $B$ be a $C^\infty$-manifold and $V$ be a smooth $\RR$-vector bundle on $B$. Let
$M\subset V$ be a vector subbundle. Then,  
the Fourier-Sato transforms on $V$ and $T_M(V)=M\oplus (V/M)$ are compatible with  specializations.
In other words, the following diagram of functors is canonically 2-commutative:
\[
\xymatrix{
D^b_\con(V) \ar[rrr]^{\FS_V} 
\ar[d]_{\nu_M}
&&& D^b(V^*)\ar[d]^{\nu_{M^\perp}}
\\
D^b(M\oplus (V/M)) \ar[rrr]_{P_{12}\circ \FS_{M\oplus (V/M)}} &&& D^b(M^\perp\oplus M^*).
}
\]
Here $P_{12}$ is the permutation of the two direct summands in $M^*\oplus M^\perp$. 

\end{thm}

The notation $\oplus$ here and below means  direct sum of vector bundles, i.e., 
fiber product over $B$. 

We note that  the diagram in Theorem \ref {refo:micro-2} can be seen as a particular
case of the outer rim of the diagram \eqref{eq:bimicro} for the case when  $X=V$,
when $M\subset V$ is our subbundle and $N=B$ is the zero section of $V$. 
In other words,  Theorem \ref {refo:micro-2} can be seen 
 as a parametrized version (with arbitrary base $B$ instead
of $B=\pt$) of 
 a particular case of Theorem \ref{thm:micro-3}
corresponding to $N=0$.

\paragraph{E. Proof of Theorem \ref{refo:micro-2}.}
The following proof is an adaptation of the argument communicated to us by P. Schapira.

\vskip .2cm

We consider three pairs
\[
 M^\perp\subset V^*, \,\,\, M\oplus M^\perp \subset V\oplus V^*, \,\,\, M\subset V, 
\]
and the corresponding normal deformations which are related by the natural projections:
\be\label{eq:three-blows}
\wt V^*_{M^\perp}\lla \wt{V\oplus V^*}_{M\oplus M^\perp} \lra \wt V_M.
\ee
Each of the three normal deformations fits into its own diagram of the form \eqref{eq:blow-diagram}
whose spaces and maps will be decorated by the subscripts $M^\perp$,  $M\oplus M^\perp$  and $M$.
In particular, the projections of the three spaces in \eqref{eq:three-blows} to the line $\RR$ will be 
denoted
  $\tau_{M^\perp}, \tau_{M\times M^\perp}$
and $\tau_M$. These projections commute with the maps in   \eqref{eq:three-blows}.
The coordinate in $\RR$ will be denoted $t$. 
 
Now, the Fourier-Sato transform on any vector bundle $W$ is defined using the region
\[
P\,\,=\,\, P_W \,\,=\,\,\bigl\{ (x,f)\in W\oplus  W^* \bigl| \,\, f(x)\geq 0\bigr\},
\]
cf. \S \ref{sec:fourier-hyp}A.  We apply this to $W=V$ and $W=M\oplus(V/M)$ and denote the
corresponding regions 
\[
P_V\subset V\oplus V^*, \quad P_{M\oplus (V/M)}\,\,\subset \,\, M\oplus  (V/M) \oplus  M^* \oplus (V/M)^*.
\]

 We want to lift $P_V$ into a region $\ol \Pc\subset \wt{V\oplus V^*}_{M\oplus M^\perp}$
 which specializes, for $t>0$, to $P_V$ and for $t=0$, to $P_{M\oplus(V/M)}$. 
 
 For this we consider the region $\Omega_{M\oplus  M^\perp}\subset \wt{V\oplus V^*}_{M\oplus M^\perp}$, 
 defined as
 the preimage $\tau_{M\oplus M^\perp}^{-1}(\RR_{>0})$, cf. \eqref{eq:blow-diagram}. It is identified with
  $V\oplus V^*\times \RR_{>0}$.
 Let $\Pc\subset\Omega$ be the image of $P_V\times \RR_{>0}$.
 
 \begin{prop}
 The closure $\ol\Pc$ of $\Pc$ in $\wt{V\oplus V^*}_{M\oplus M^\perp}$ is the union of $\Pc$
 and $P_{M\oplus (V/M)}\subset \tau_{M\oplus M^\perp}^{-1}(0)$. 
 \end{prop}
 
 \noindent{\sl Proof:}  The statement is local in $B$. So we can assume that there exists  
  a complement $M'$ to $M$ and to write $V=M\oplus M'$.  We then
  identify, as in \eqref{eq:blow-diagram},
 \[
 \wt{V\oplus V^*}_{M\oplus M^\perp} \,\,=\,\, M\oplus M'\oplus M^* \oplus  M'{}^*\times\RR
 \]
 and the projection $p_{M\oplus M^\perp}:  \wt{V\oplus  V^*}_{M\oplus M^\perp}\to V\oplus  V^*$ can be written as
 \be\label{eq:p-M-perp}
 \begin{gathered}
 p_{M\oplus M^\perp}: M\oplus  M'\oplus  M^* \oplus  M'{}^*\times\RR \lra M\oplus M'\oplus M^* \oplus M'{}^*, 
 \\
 \bigl( m,m', \phi, \phi', t \bigr) \,\,\mapsto \,\, \bigl(m, tm', t\phi, \phi'\bigr). 
 \end{gathered}
 \ee
 Recall that the identification $\Omega_{M\oplus  M^\perp}\to V\oplus V^*\times \RR_{>0}$
 is given by the map $(p_{M\oplus M^\perp}, \tau_{M\oplus  M^\perp})$, the second component being
 projection to $t$. It follows from \eqref{eq:p-M-perp} that
 for any $t>0$ the image, under $p_{M\oplus  M^\perp}$,  of $P_V\times\{t\}$ is $P_V$. 
 Therefore 
  the inverse of 
 $p_{M\oplus  M^\perp}, \tau_{M\oplus  M^\perp})$
 identifies $P_V\times \RR_{\geq 0}$ with $P_V\times \RR_{\geq 0}$, where for $t=0$ our choice of
 complement has identified $P_V$ with $P_{M\oplus  M'} = P_{M\oplus (V/M)}$. \qed

 \vskip .2cm
 
 We now consider the following diagram:
 \[
 \xymatrix{
 V^* && \ar[ll]_{p_{2,V}} P_V \ar[rr]^{p_{1,V}}&& V
 \\
 \Omega_{M^\perp} 
 \ar[u]^{\wt p_{M^\perp}} 
 \ar[d]_{j_{M^\perp}}
 && 
 \ar[ll]_{\rho_2} \Pc
 \ar[u]
  \ar[rr]^{\rho_1}
  \ar[d]
  && \Omega_M
 \ar[u]_{\wt p_M}
 \ar[d]^{j_M}
 \\
 \wt V^*_{M^\perp}&& \ar[ll]_{\pi_2} \wt\Pc \ar[rr]^{\pi_1} && \wt V_M
 \\
 M^\perp\oplus M^*
 \ar[u]^{s_{M^\perp}} &&\ar[ll]^{p_{2, M\oplus (V/M)}} P_{M\oplus (V/M)} \ar[rr]_{p_{1, M\oplus (V/M)}} 
 \ar[u]
 &&
 M\oplus (V/M) \ar[u]_{s_M}
 }
 \]
 Given $\Gc\in D^b_\con(V)$,  we have that
 \[
 \begin{gathered}
 \nu_{M^\perp}\FS_V(\Gc) \,\,=\,\, s_{M^\perp}^* R(j_{M^\perp})_* \,\wt p_{M^\perp}^* \,(p_{2,V})_! \, p_{1,V}^{-1} (\Gc),
 \\
 \FS_{M\oplus (V/M)} \nu_M(\Gc) \,\,=\,\,    (p_{2, M\oplus  (V/M)}) _! \, \, 
 p_{1, M\oplus  (V/M)}^*\, 
    s_M^* \, R(j_M)_* \,  \wt p_M^* (\Gc)
 \end{gathered}
 \]
 are given by moving along the two boundary paths of this diagram from the northeast to the southwest corner. 
 We identify these functors using the base change theorem for the Cartesian squares forming this diagram.

\paragraph{F. Proof of Theorem \ref{thm:micro-3}.}
We write the diagram \eqref {eq:bimicro} in our case as follows:

\be\label{eq:MMN}
\xymatrix{
D^b(V, \Sc_\RR) \ar[rr]^{\mu_M} 
\ar[d]_{\mu_N}
\ar[rrd]^{\mu_{NM}}
&& D^b_\con(M \times (V/M)^*)
\ar[d]^{\mu_{N\times(V/M)^*}}
\\
D^b_\con(N\times (V/N)^*) \ar[rr]_{\hskip -1cm \nu_{N\times (V/M)^*}}&& 
D^b_\bic(N\times (M/N)^*\times (V/M)^*).
}
\ee
Here  and below the subscript ``con'' means complexes which are $\RR_{>0}$-conic with to the second
argument, and ``bico'' means cmplexes which are $(\RR_{>0})^2$-biconic with respect to the
second and third arguments.

We recall that $\mu_{NM}$ is the composition
\[
D^b(V, \Sc_\RR) \buildrel \nu_{NM}\over\lra D^b_\bic(N\times (M/N) \times (V/M))
\buildrel \FS_{(M/N)\times (V/M)}\over\lra D^b_\bic(N\times (M/N)^* \times (V/M)^*).
\]
We now prove the 2-commutativity of each of the two triangles in \eqref{eq:MMN}. 

\vskip .2cm

\noindent \underbar {\sl Upper triangle.} We write each $\mu$ as the composition of the corresponding $\FS$ and $\nu$
and apply Theorem \ref{thm:bispec} to decompose $\nu_{NM}$ as the composition of two specializations. After this
we represent the two paths in the triangle as the two boundary paths in the following diagram:
\[
\xymatrix{
D^b(V, \Sc_\RR) 
\ar[d]_{\nu_M} 
\ar[dr]^{\mu_M}
&
\\
D^b_\con(M\times(V/M)) \ar[r]^{\FS_{V/M}}
\ar[d]_{\nu_N}
& D^b_\con(M\times (V/M)^*)
\ar[d]^{\nu_N}
\\
D^b_\bic(N\times(M/N)\times(V/M)) \ar[r]^{\FS_{V/M}}
\ar[dr]_{\FS_{(M/N)\times (V/M)}}
&
 D^b_\bic(N\times (M/N) \times(V/M)^*))
 \ar[d]^{\FS_{M/N}}
 \\
 &D^b_\bic(N\times(M/N)^* \times (V/M)^*). 
}
\]
In this diagram, the top triangle commutes by definition of $\mu_M$ and  the commutativity of the  bottom triangle
expresses the fact that the Fourier-Sato transform  of biconic sheaves on the  direct sum of vector bundles
 can be done in stages, cf. \cite{KaSha} Prop. 3.7.15. The commutativity 
of the middle square follows because specialization along $N$ and the Fourier-Sato transform along $V/M$
operate in different factors so they are independent of each other and  can be permuted. 

\vskip .2cm

\noindent\underbar{\sl Lower triangle.} As before, by unravelling the definitions of various $\mu$ and
applying Theorem  \ref{thm:bispec}, we represent the two paths in the triangle as the two boundary paths in the following diagram:
\[
\xymatrix{
&& D^b(V, \Sc_\RR) 
\ar[d]^{\nu_N}
\ar[lld]_{\mu_N}
\\
D^b_\con(N\times(V/N)^*)
\ar[d]_{\nu_{N\times(V/M)^*}}
&&  \ar[ll]_{\FS_{V/N}}D^b_\con(N\times(V/N))
\ar[d]^{\nu_{N\times(M/N)}}
\\
D^b_\bic(N\times(M/N)\times(V/M)) && \ar[ll]^{\FS_{(M/N)\times (V/M)}} D^b_\bic(N\times (M/N) \times (V/M))
}
\]
The commutativity of the top triangle in this diagram is the definition of $\mu_N$. The commutativity of
the lower square is an instance of Theorem \ref{refo:micro-2} for the  trivial vector bundle over $B=N$ with
fiber  $V/N$ and the trivial subbundle with fiber $M/N$. 
Theorem \ref{thm:micro-3} is proved.

    \vfill\eject

  \vskip .5cm
  
  {\small
  
    M.F.: Department of Mathematics, National Research University Higher School
    of Economics, Russian Federation, 6 Usacheva st, Moscow 119048;
    Skolkovo Institute of Science and Technology;
    Institute for Information Transmission Problems;
{\tt fnklberg@gmail.com}  
  
    M.K.: Kavli IPMU, 5-1-5 Kashiwanoha, Kashiwa, Chiba, 277-8583 Japan,
{\tt mikhail.kapranov@ipmu.jp}

V.S.: Institut de Math\'ematiques de Toulouse, Universit\'e Paul Sabatier, 118 route de Narbonne, 
31062 Toulouse, France, 
 {\tt schechtman@math.ups-tlse.fr }
  
 }

\end{document}